\documentclass[a4paper, 11pt]{amsart}
\usepackage[utf8]{inputenc}

\usepackage[english]{babel}

\usepackage{fancyheadings}
\usepackage{amsmath,amssymb,amsthm}
\usepackage{mathrsfs}
\usepackage{tikz}

\usepackage{a4wide}
\usepackage[text={16cm,24cm}]{geometry}
\usepackage{graphicx}
\usepackage[arrow, matrix, curve]{xy}
\usepackage{MnSymbol}

\numberwithin{equation}{section}

\theoremstyle{definition}
\newtheorem{df}{Definition}[section]
\newtheorem{rem}[df]{Remark}
\newtheorem{ex}[df]{Example}

\theoremstyle{plain}
\newtheorem{lemma}[df]{Lemma}
\newtheorem{prop}[df]{Proposition}
\newtheorem{thm}[df]{Theorem}
\newtheorem{cor}[df]{Corollary}

\newenvironment{bew}{\begin{proof}[Proof]}{\end{proof}}

\renewcommand{\a}{\mathbb A}
\renewcommand{\c}{\mathbb C}
\newcommand{\f}{\mathbb F}
\newcommand{\ka}{\mathbb K}

\newcommand{\p}{\mathbb P}

\renewcommand{\t}{\mathbb T}
\newcommand{\z}{\mathbb Z}

\newcommand{\F}{\mathcal F}
\newcommand{\G}{\mathcal G}

\newcommand{\M}{\mathfrak M}

\newcommand{\Jon}{\text{\textbf{Jon}}}
\newcommand{\Aff}{\text{\textbf{Aff}}}

\DeclareMathOperator{\Aut}{Aut}
\DeclareMathOperator{\SAut}{SAut}
\DeclareMathOperator{\supp}{\text{supp}}

\DeclareMathOperator{\ML}{\text{ML}}

\DeclareRobustCommand{\ie}{i.\,e.~}

\newcommand{\nlin}{\unitlength1mm\begin{picture}(0,9.25)
                       \put(0,0.75){\line(0,1){8.5}}
                      \end{picture}}

\newcommand{\vlin}[1]{\hspace{0.75mm}\unitlength1mm\begin{picture}(#1,0)
                       \put(0,0){\line(1,0){#1}}
                      \end{picture}\hspace{0.75mm}\rule[-3mm]{0mm}{4mm}}

\newcommand{\lin}{\vlin{8.5}}

\newcommand{\lllin}{\vlin{15}}

\newcommand{\co}[1]{\unitlength1mm\begin{picture}(0,8)
    \put(0,0){\circle{1.5}}
    \put(0,3){\makebox(0,5)[b]{$#1$}}
                      \end{picture}}

\newcommand{\mybox}{\unitlength1mm\begin{picture}(0,1.5)
    \put(-0.75,-0.75){\line(0,1){1.5}}
    \put(-0.75,-0.75){\line(1,0){1.5}}
    \put(0.75,0.75){\line(0,-1){1.5}}
    \put(0.75,0.75){\line(-1,0){1.5}}
    \end{picture}}
    
\newcommand{\boxo}[1]{\unitlength1mm\begin{picture}(0,8)
    \put(0,0){\mybox}
    \put(0,3){\makebox(0,5)[b]{$#1$}}
                      \end{picture}}

\newcommand{\xbox}{\unitlength1mm\begin{picture}(0,1.5)
    \put(0,0){$\mybox$}
    \put(-0.75,0){\line(1,0){1.5}}
    \put(0,-0.75){\line(0,1){1.5}}
    \end{picture}}
    
\newcommand{\xboxo}[1]{\unitlength1mm\begin{picture}(0,8)
    \put(0,0){\xbox}
    \put(0,3){\makebox(0,5)[b]{$#1$}}
                      \end{picture}}

\newcommand{\cu}[1]{\unitlength1mm\begin{picture}(0,8)
    \put(0,0){\circle{1.5}}
    \put(0,-7){\makebox(0,4)[b]{$#1$}}
    \end{picture}
      \rule[-7mm]{0mm}{7mm}}

\newcommand{\cou}[2]{\unitlength1mm\begin{picture}(0,8)
    \put(0,0){\circle{1.5}}
    \put(0,3){\makebox(0,5)[b]{$#1$}}
    \put(0,-8.5){\makebox(0,4)[t]{$#2$}}
      \end{picture}
      \rule[-7mm]{0mm}{7mm}}

\newcommand{\cshiftup}[2]{\unitlength1mm\begin{picture}(0,9.25)
                       \put(0,10){\cou{#1}{#2}}
                      \end{picture}}

\newcommand{\xbou}[2]{\unitlength1mm\begin{picture}(0,8)
    \put(0,0){\xbox}
    \put(0,2){\makebox(0,5)[b]{$#1$}}
   \put(0,-7){\makebox(0,5)[b]{$#2$}}
      \end{picture}
      \rule[-7mm]{0mm}{7mm}}

\newcommand{\xbshiftup}[2]{\unitlength1mm\begin{picture}(0,9.25)
                       \put(0,10){\xbou{#1}{#2}}
                      \end{picture}}

\begin{document}

\begin{center}
\LARGE{Transitivity of automorphism groups of Gizatullin surfaces}

\parskip 20pt

\normalsize Sergei Kovalenko
\end{center}

\vspace{10pt}

\begin{sloppypar}
\noindent\small ABSTRACT. We show that the automorphism group of a certain subclass of smooth Gizatullin surfaces with a distinguished and rigid extended divisor is generated by automorphisms of $\a^1$-fibrations. Moreover, such surfaces yield examples of smooth Gizatullin surfaces with a non-transitive action of the automorphism group. Thus, they represent counter-examples to Gizatullin's conjecture. For such surfaces we give an explicit orbit decomposition of the natural action of the automorphism group in some special cases. 
\end{sloppypar}

\tableofcontents

\begin{sloppypar} 

\section{Introduction}

In the following we consider surfaces over the field $\ka = \c$ of complex numbers. $\a^n$ stands for the affine $n$-space over $\c$. All results also hold for arbitrary algebraically closed fields of characteristic zero.\\

Gizatullin surfaces were introduced by Danilov and Gizatullin (\cite{Gi}, \cite{DG2} and \cite{DG3}). We recall that the Makar-Limanov invariant $\ML(V)$ of a normal affine surface $V$ is defined by

$$\ML(V) = \bigcap_{\partial \ \in \ \text{LND}(\c[V])} \ker(\partial).$$ 

A useful characterization of normal affine surfaces with trivial Makar-Limanov invariant is the following result due to Gizatullin (\cite{Gi}, Theorems 2 and 3), Bertin (\cite{Be}, Theorem 1.8), Bandman and Makar-Limanov (\cite{BML} and \cite{ML}) in the smooth case and due to Dubouloz (\cite{Du}) in the normal case:

\begin{prop}\label{CharacterizationGizatullinSurface}(\cite{FZ}, Theorem 4.3) For a normal affine surface $V$ that is non-isomorphic to $\c^* \times \c^*$ or to $\c^* \times \a^1$, the following conditions are equivalent:
\begin{itemize}
\item[(1)] $\ML(V)$ is trivial, \ie $\ML(V) = \c$.
\item[(2)] The automorphism group acts on $V$ with an open orbit, such that the complement is finite (such orbits are called \emph{big}).
\item[(3)] $V$ admits a smooth completion by a zigzag $D$. In other words, $V = X \backslash D$, where $X$ is a complete surface smooth along $D$ and $D = C_0 \cup \cdots \cup C_n$ is a linear chain of smooth rational curves with simple normal crossings.
\end{itemize}
\end{prop}

Normal affine surfaces $V$ satisfying one of the equivalent conditions of Proposition \ref{CharacterizationGizatullinSurface} are called \emph{Gizatullin surfaces}. In particular, the automorphism group of Gizatullin surface $V$ is quite large compared to surfaces in general. Studying simple smooth Gizatullin surfaces, like Danielewski surfaces $V_P = \{ xy - P(z) = 0 \} \subseteq \a^3$ with a polynomial $P$ having pairwise distinct roots, one easily sees that the automorphism group acts transitively on these surfaces (see \cite{ML}). In this case, the big orbit $O$ coincides with $V$. Cleary, in the singular case $O$ cannot coincide with $V$. However, it is still an open question whether $O$ coincides with $V$ in the smooth case. More generally Gizatullin formulated the following conjecture in \cite{Gi}:\\

\noindent \textbf{Conjecture (\cite{Gi}, Conjecture 1)\footnote[1]{See Remark \ref{OrevkovCurves}.}:} Let $V$ be a smooth Gizatullin surface over an algebraically closed field $\ka$ of characteristic zero. Then the action of the automorphism group is transitive on $V$. In other words, the open orbit $O$ of $\Aut(V)$ coincides with $V$.\\

In contrast to the case of characteristic zero, counter-examples in positive characteristic were found very early and can be found in \cite{DG1}.\\

In general it is difficult to determine the orbits of the natural action of $\Aut(V)$ for general normal affine varieties $V$. But there are some nice results in higher dimension. For example, in \cite{DM-JP} it is shown that the automorphism group of the Koras-Russell cubic $X = \{ x + x^2y + z^2 + t^3 = 0 \} \subseteq \a^4$ has exactly $4$ orbits, one of them being the fixed point $p = (0, 0, 0, 0)$.\\

The aim of this article is to construct families of smooth Gizatullin surfaces $V$, satisfying special conditions and to determine the orbit decomposition of the natural action of $\Aut(V)$. In particular, it turns out that for these surfaces the big orbit $O$ is a proper subset of $V$. It follows that such Gizatullin surfaces provide counter-examples to the Gizatullin conjecture.

Let $(X, D)$ be an SNC-completion of a Gizatullin surface $V$ so that $V = X \backslash D$ and $D$ is a simple normal crossing divisor. It is well known that $D$ can be transformed by birational transformations into \emph{standard form}. This means that $D = C_0 \cup \cdots \cup C_n$ is a chain of smooth rational curves satisfying $C_0^2 = C_1^2 = 0$ and $C_i^2 \leq -2$ for $i \geq 2$ if $n \geq 4$ or $C_i^2 = 0$ for all $i$ if $n \leq 3$. The standard form of the boundary divisor $D$ is (up to reversion) an invariant of the abstract isomorphism type of $V$ (see \cite{FKZ1C}, Cor. 3.33'). But in general this invariant is too weak. It is more convenient to introduce a stronger invariant, the so called \emph{extended divisor} $D_{\text{ext}}$, defined as follows. The linear pencils $|C_0|$ and $|C_1|$ provide $\p^1$-fibrations $\Phi_0 := \Phi_{|C_0|}: X \to \p^1$ and $\Phi_1 := \Phi_{|C_1|}: X \to \p^1$. These $\p^1$-fibrations lift to the minimal resolution of singularities $\tilde{X}$ of $X$. By [FKZ2], Lemma 2.19, $\Phi_0$ admits at most one degenerate fiber, without lost of generality the fiber over $0$, and the \emph{extended divisor} of $(X, D)$ is

$$D_{\text{ext}} := C_0 \cup C_1 \cup \Phi_0^{-1}(0).$$

The extended divisor $D_{\text{ext}}$ always contains the boundary divisor $D$. The connected components of $D_{\text{ext}} - D$ are called \emph{feathers}. We denote them by $F_{i, j}$, $2 \leq i \leq n$, $j \in \{ 1, \dots, r_i \}$, and assume that $F_{i, j}$ is attached to the curve $C_i$ at points $P_{i, j}$. In fact, in the smooth case the feathers are irreducible. The Matching Principle (see \cite{FKZ4}) states that there is a natural bijection between feathers $F_{i, j}$ of $(X, D)$ and feathers $F^\vee_{i, j}$ of the completion $(X^\vee, D^\vee)$ obtained by reversing the boundary zigzag. Our candidates for potential counter-examples for the Gizatullin conjecture are smooth Gizatullin surfaces $V$ with the property that there exist a standard completion $(X, D)$, $D = C_0 \cup \cdots \cup C_n$ of $V$ such that $C_3, \dots, C_{n - 1}$ are \emph{inner components} (see Def. \ref{StarComponent}) and such that no feathers are attached to $C_2$ and to $C_n$. For a precise formulation of the result see Theorem \ref{FixedPointOfSurface} in section 3. 

Furthermore, we show that in the general case the automorphism group of such surfaces $V$ is generated by automorphisms of $\a^1$-fibrations, that is, by automorphisms which preserve certain $\a^1$-fibrations. We also give an explicit description of $\Aut(V)$ as an amalgamated product of two subgroups.

\bigskip
This article is structured as follows. In section 2 we introduce the main tools employed to work with Gizatullin surfaces and more generally with $\a^1$-fibered surfaces. We are mainly interested in presentations and properties of standard and $1$-standard completions of such surfaces. In particular, we give a decomposition of birational maps between $1$-standard pairs. Following \cite{FKZ3}-\cite{FKZ4}, we introduce the Matching Principle and the rigidity of extended divisors. Finally, we give a concrete description of smooth Gizatullin surfaces in local affine coordinates.

In section 3 we apply these tools to show that the automorphism group of a general smooth Gizatullin surface $V$ with a distinguished and rigid extended divisor is generated by automorphisms of $\a^1$-fibrations. Moreover, for such surfaces we give the orbit decomposition of the action of the automorphism group of $V$. In particular we show that the automorphism group admits fix points in general. These surfaces provide counter-examples to the Gizatullin conjecture. In subsection \ref{AmalgamatedProductStructure} we give explicit presentations of the automorphism groups of such surfaces as amalgamated products of two automorphism subgroups.

Finally, in Section 4 we deal with the singular case and we show that similar results hold.

\bigskip
\noindent \textbf{Acknowledgements} The results presented in this paper are contained in the author's Ph.D. thesis \cite{Ko}. The author would like to thank his supervisor Professor Dr. Hubert Flenner for inspiring discussions and mathematical advice. Furthermore, the author thanks S. Kaliman for detailed suggestions and remarks for fixing the proof of the main result. The author also thanks Anne Wald and Felix Fr\"uhauf for helpful comments and remarks.

\end{sloppypar}

\begin{sloppypar} 

\section{Preliminaries}

\subsection{$\a^1$-fibered surfaces and Gizatullin surfaces}

Following \cite{DG2}, we introduce the notion of an oriented zigzag:

\begin{df} A \emph{zigzag} $D$ on a normal projective surface $X$ is an SNC-divisor supported in the smooth locus $X_{\text{reg}}$ of $X$, with irreducible components isomorphic to $\p^1$ and whose dual graph is a chain. If $\supp(D) = \bigcup_{i = 0}^n C_i$ is the decomposition into irreducible components, one can order the $C_i$ such that

$$C_i.C_j = \begin{cases} 1, & |i - j| = 1 \\ 0, & |i - j| > 1 .\end{cases}$$

\noindent A zigzag with such an ordering is called \emph{oriented} and the sequence $[[(C_0)^2, \dots, (C_n)^2]]$ is called the \emph{type} of $D$. The same zigzag with the reverse ordering is denoted by ${}^tD$ (\ie ${}^tD$ is of type $[[(C_n)^2, \dots, (C_0)^2]]$).\\
An \emph{oriented sub-zigzag} of an oriented zigzag is an SNC-divisor $D'$ with $\supp(D') \subseteq \supp(D)$ which is a zigzag for the induced ordering.\\
We say that an oriented zigzag $D$ is composed of sub-zigzags $Z_1, \dots, Z_s$, and following [BD] we denote $D = Z_1 \triangleright \cdots \triangleright Z_s$, if the $Z_i$, $1\leq i\leq s$, are oriented sub-zigzags of $D$ whose union is $D$ and the components of $Z_i$ precede those of $Z_j$ for $i < j$.
\end{df}

Surfaces completable by a zigzag were first studied by Danilov and Gizatullin (\cite{Gi}, \cite{DG2} and \cite{DG3}).

\begin{df} A normal affine surface $V$ is called a \emph{Gizatullin surface} if it is completable by a zigzag.
\end{df}

For the rest of this article we fix the following notation:

\bigskip
\noindent \textbf{Notation:} If $V$ is a Gizatullin surface and $(X, D)$ is a completion of $V$ by a zigzag $D$, then

$$D = C_0 + \cdots + C_n \quad \text{and} \quad C_i \ \text{and} \ C_j \ \text{have a non-empty intersection only for} \ |i - j| = 1.$$

\bigskip
\noindent Given a Gizatullin surface $V$ together with a completion $(X, D)$ by a zigzag, we can associate a linear weighted graph $\Gamma_D$ to $(X, D)$. The vertices $v_i$, $0 \leq i \leq n$, are the boundary components $C_i$ and the weights are the corresponding self-intersection numbers $w_i := C_i^2$. Thus $\Gamma_D$ has the form

$$\Gamma_D: \quad \cou{C_0}{w_0} \lin \cou{C_1}{w_1} \lin \ \cdots \ \lin \cou{C_n}{w_n} \quad .$$

\noindent For a better systematic understanding of Gizatullin surfaces we introduce elementary transformations of weighted graphs.

\begin{df} Given an at most linear vertex $v$ of a weighted graph $\Gamma$ with weight $0$ one can perform the following transformations. If $v$ is linear with neighbors $v_1, v_2$ then we blow up the edge connecting $v$ and $v_1$ in $\Gamma$ and blow down the proper transform of $v$:

\begin{equation}\label{ElementaryTransformation1}
\dots \ \cou{v_1}{w_1 - 1} \lin \cou{v'}{0} \lin \cou{v_2}{w_2 + 1} \ \dots \ \dasharrow \ \dots \ \cou{v_1}{w_1 - 1} \lin \cou{v'}{-1} \lin \cou{v}{-1} \lin \cou{v_2}{w_2} \ \dots \ \to \ \dots \ \cou{v_1}{w_1} \lin \cou{v}{0} \lin \cou{v_2}{w_2} \quad .
\end{equation}

\noindent Similarly, if $v$ is an end vertex of $\Gamma$ connected to the vertex $v_1$ then one proceeds as follows:

\begin{equation}\label{ElementaryTransformation2}
\dots \ \cou{v_1}{w_1 - 1} \lin \cou{v'}{0} \ \dasharrow \ \dots \ \cou{v_1}{w_1 - 1} \lin \cou{v'}{-1} \lin \cou{v}{-1} \ \to \ \dots \ \cou{v_1}{w_1} \lin \cou{v}{0} \quad .
\end{equation}

\noindent These operations (\ref{ElementaryTransformation1}) and (\ref{ElementaryTransformation2}) and their inverses are called \emph{elementary transformations} of $\Gamma$. If such an elementary transformation involves only an inner blowup then we call it \emph{inner}. Thus (\ref{ElementaryTransformation1}) and (\ref{ElementaryTransformation2}) are inner whereas the inverse of (\ref{ElementaryTransformation2}) is not as it involves an outer blowup.
\end{df}

\bigskip
We consider a Gizatullin surface $V = X \backslash D$, where $X$ is projective and $D$ is a zigzag. By a sequence of blowups and blowdowns we can transform the dual graph $\Gamma_D$ of $D$ into \textit{standard form}, \ie we can achieve that $C_0^2 = C_1^2 = 0$ and $C_i^2 \leq -2$ for all $i \geq 2$ if $n \geq 4$ or $C_i^2 = 0$ for all $i$ if $n \leq 3$ (cf. \cite{DG2}, \cite{Da}, \cite{FKZ1}). Moreover, this representation is unique up to reversion meaning that for two standard forms $[[0, 0, w_2, \dots, w_n]]$ and $[[0, 0, w'_2, \dots, w'_n]]$ either $w_i = w'_i$ or $w_i = w'_{n + 2 - i}$ holds (\cite{FKZ1C}).

The reversion process can be described as follows. We start with a boundary divisor of type $[[0, 0, w_2, \dots, w_n]]$. Performing the elementary transformation (\ref{ElementaryTransformation1}) at the vertex corresponding to $C_1$ we obtain a boundary divisor of type $[[-1, 0, w_2 + 1, w_3, \dots, w_n]]$. After $|w_2|$ steps we arrive at a boundary divisor of type $[[w_2, 0, 0, w_3, \dots, w_n]]$. This means that we can move pairs of zeros to the right. Repeating this, we finally obtain a boundary divisor of type $[[w_2, \dots, w_n, 0, 0]]$. Notice that all birational transformations are centered in the boundary, \ie these transformations yield the identity on the affine parts.\\

We recall the notion of an $m$-standard zigzag (see \cite{DG2}, (1.2)):

\begin{df} A zigzag $D$ on a normal projective surface $X$ is called \emph{$m$-standard} (or in \emph{$m$-standard form}), if it is of type $[[0, -m, w_2, \dots, w_n]]$ with $n \geq 1$ and $w_i \leq -2$ (in the case of $n = 1$ there are no weights $w_i$).\\
An \emph{$m$-standard pair} is a pair $(X, D)$ consisting of a normal projective surface $X$ and an $m$-standard zigzag $D$ on $X$. If $m = 0$, then $(X, D)$ is called a \emph{standard pair}. A \emph{birational map} $\varphi: (X, D) \dasharrow (X', D')$ between $m$-standard pairs is a birational map $\varphi: X \dasharrow X'$ which restricts to an isomorphism $\varphi\vert_{X \backslash D}: X \backslash D \stackrel{\sim}{\to} X' \backslash D'$.
\end{df}

Since the underlying projective surface $X$ of an $m$-standard pair is rational, it is equipped with a rational fibration $\bar{\pi} = \Phi_{|C_0|}: X \to \p^1$ defined by the complete linear system $|C_0|$. In particular, if $m = 0$, there are even two $\p^1$-fibrations $\Phi_0 := \Phi_{|C_0|}, \Phi_1 := \Phi_{|C_1|} : X \to \p^1$ and, thus, a morphism 

$$\Phi := \Phi_0 \times \Phi_1: X \to \p^1 \times \p^1,$$

\noindent which is birational (\cite{FKZ2}, Lemma 2.19). After a change of coordinates we can assume that $C_0 = \Phi_0^{-1}(\infty)$, $\Phi(C_1) = \p^1 \times \{ \infty \}$ and $C_2 \cup \cdots \cup C_n \subseteq \Phi_0^{-1}(0)$. The divisor $D_{\text{ext}} := C_0 \cup C_1 \cup \Phi_0^{-1}(0)$ is called the \textit{extended divisor}. We also denote the full fiber $\Phi_0^{-1}(0)$ by $D_{(e)}$. Before determining the structure of the extended divisor, we recall the notion of a \emph{feather}:

\begin{df} (\cite{FKZ2}, Def. 5.5)
\begin{itemize}
\item[(1)] A \emph{feather} is a linear chain 

$$F: \quad \co{B}{} \lin \co{F_1} \lin \dots \lin \co{F_s}$$

\noindent of smooth rational curves such that $B^2 \leq -1$ and $F_i^2 \leq -2$ for all $i \geq 1$. The curve $B$ is called the \emph{bridge curve}.
\item[(2)] A collection of feathers $\{ F_\rho \}$ consists of feathers $F_\rho$, $1 \leq \rho \leq r$, which are pairwise disjoint. Such a collection will be denoted by a plus box 

$$\xboxo{ \{ F_\rho \} } \quad .$$

\item[(3)] Let $D = C_0 + \cdots + C_n$ be a zigzag. A collection $\{ F_\rho \}$ is \emph{attached to a curve} $C_i$ if the bridge curves $B_\rho$ meet $C_i$ in pairwise distinct points and all the feathers $F_\rho$ are disjoint with the curves $C_j$ for $j \neq i$.
\end{itemize} 
\end{df}

\begin{lemma}\label{StructureExtendedDivisor} (\cite{FKZ3}, Prop. 1.11) Let $(\tilde{X}, D)$ be a minimal SNC completion of the minimal resolution of singularities of a Gizatullin surface $V$. Furthermore, let $D = C_0 + \cdots + C_n$ be the boundary divisor in standard form. Then the extended divisor $D_{\text{ext}}$ has the dual graph

\vspace{15pt}
$$D_{\text{ext}}: \quad \cou{0}{C_0} \lin \cou{0}{C_1} \lin \cu{C_2} \nlin \xbshiftup{ \{ F_{2, j} \} }{} \lin \dots \lin \cu{C_i} \nlin \xbshiftup{ \{ F_{i, j} \} }{} \lin \dots \lin \cu{C_n} \nlin \xbshiftup{ \{ F_{n, j} \} }{} \quad ,$$

\noindent where $\{ F_{i, j} \}$, $j \in \{ 1, \dots, r_i \}$, are feathers attached to the curve $C_i$. Moreover, $\tilde{X}$ is obtained from $\p^1 \times \p^1$ by a sequence of blowups with centers in the images of the components $C_i$, $i \geq 2$.
\end{lemma}

\begin{rem} We consider the feathers $F_{i, j} := B_{i, j} + F_{i, j, 1} + \cdots + F_{i, j, k_{i, j}}$ mentioned in Lemma \ref{StructureExtendedDivisor}. The collection of linear chains $R_{i, j} := F_{i, j, 1} + \cdots + F_{i, j, k_{i, j}}$ corresponds to the minimal resolution of singularities of $V$. Thus, if $(X, D)$ is a standard completion of $V$ and $(\tilde{X}, D)$ is the minimal resolution of singularities of $(X, D)$, the chain $R_{i, j}$ contracts via $\mu: (\tilde{X}, D) \to (X, D)$ to a singular point of $V$, which is a cyclic quotient singularity. In partcular, $V$ has at most cyclic quotient singularities ([Mi], \S 3, Lemma 1.4.4 (1) and \cite{FKZ3}, Remark 1.12).

Hence, $V$ is smooth if and only if every $R_{i, j}$ is empty, \ie if every feather $F_{i, j}$ is irreducible and reduces to a single bridge curve $B_{i, j}$ ([FKZ3], 1.8, 1.9 and Remark 1.12).
\end{rem}

In connection with Lemma \ref{StructureExtendedDivisor} we abbreviate the subdivisor $\sum_{k \geq i} C_k + \sum_{j_k; k \geq i} F_{k, j_k}$ by $D^{\geq i}_{\text{ext}}$ and the subdivisor $\sum_{k > i} C_k + \sum_{j_k; k \geq i} F_{k, j_k} = D^{\geq i}_{\text{ext}} \ominus C_i$ by $D^{> i}_{\text{ext}}$.

\bigskip
Similarly as standard completions of Gizatullin surfaces arise from $\p^1 \times \p^1$, $1$-standard completions arise from the Hirzebruch surface $\f_1$. More explicitly, we have the following lemma:

\begin{lemma}\label{FactorizationHirzebruch}(\cite{BD}, Lemma 1.0.7) Let $(X, D)$ be a $1$-standard pair and let $\mu: \tilde{X} \to X$ be the minimal resolution of singularities of $X$. Then there exists a birational morphism $\eta: \tilde{X} \to \f_1$, unique up to an automorphism of $\f_1$, that restricts to an isomorphism outside the degenerate fibers of $\bar{\pi} \circ \mu$, and satisfies the commutative diagram

$$\begin{xy}
  \xymatrix{
   & \tilde{X} \ar[ld]_\mu \ar[rd]^\eta \ar[dd]^{\mu \circ \bar{\pi}} & \\
  X \ar[rd]_{\bar{\pi}} &  & \f_1 \ar[ld]^\rho \\
   & \p^1 & .
  }
\end{xy}$$

\noindent Moreover, if $(X', D')$ is another $1$-standard pair with associated morphism $\eta': \tilde{X}' \to \f_1$, then $(X, D)$ and $(X', D')$ are isomorphic if and only if there exists an automorphism of $\f_1$ isomorphically mapping $\eta(\mu^{-1}_*(C_0))$ onto $\eta'(\mu'^{-1}_*(C'_0))$ and isomorphically sending the base-points of $\eta^{-1}$ (including infinitely near ones) onto those of $\eta'^{-1}$.
\end{lemma}

\bigskip
Since $\f_1$ is the blowup of $\p^2$ in one point, every $1$-standard pair $(X, D)$ arises as a blowup of $\p^2$ and the blowup process starts as follows

$$\cdots \quad \rightarrow \quad \co{0} \lin \co{-1} \lin \co{0} \quad \rightarrow \quad \co{1} \lin \co{1} \quad .$$

\noindent Here we can take any two lines in $\p^2$ for the two curves with self-intersection $1$. The extended divisor can also be defined for $1$-standard pairs and $D_{\text{ext}}$ becomes\\

\bigskip
$$D_{\text{ext}}: \quad \cou{0}{C_0} \lin \cou{-1}{C_1} \lin \cu{C_2} \nlin \xbshiftup{ \{ F_{2,j_2} \} }{} \lin \dots \lin \cu{C_i} \nlin \xbshiftup{ \{ F_{i,j_i} \} }{} \lin \dots \lin \cu{C_n} \nlin \xbshiftup{ \{ F_{n,j_n} \} }{} \quad .$$

\noindent This results in the same divisor as taking the extended divisor of the corresponding standard completion, blowing up the intersection point $C_0 \cap C_1$ and blowing down the proper transform of $C_0$.
 
We will often deal with $1$-standard pairs. It follows from \cite{BD}, Lemma 2.1.1 that every birational map $\varphi: (X, D) \dasharrow (X', D')$ between $1$-standard pairs which is not an isomorphism has a unique base point $p \in C_0$. This base point is called the \textit{center} of $\varphi$. In general, this yields qualitatively different maps depending on whether $p \in C_0 \cap C_1$ or $p \in C_0 \backslash C_1$.

\begin{df} Two $\a^1$-fibered surfaces $(V, \pi)$ and $(V', \pi')$ are said to be isomorphic if there exists an isomorphism $\Psi: V \to V'$ and an automorphism $\psi$ of $\a^1$, such that $\pi' \circ \Psi = \psi \circ \pi$.\\
Two $\a^1$-fibrations $\pi, \pi'$ on a surface $V$ are said to be isomorphic if $(V, \pi)$ and $(V, \pi')$ are isomorphic.
\end{df}

As mentioned above, there are two basic types of birational maps between $1$-standard pairs: The \textit{fibered modifications}, which preserve the given $\a^1$-fibrations, and the \textit{reversions}, which are, in some sense, the simplest maps that do not preserve the given fibrations.

\begin{df} Let $\varphi: (X, D) \dasharrow (X', D')$ be a birational map between $1$-standard pairs and let $D = C_0 \triangleright \cdots \triangleright C_n$ and $D' = C'_0 \triangleright \cdots \triangleright C'_n$ be the oriented boundary divisors.
\begin{itemize}
\item[(1)] (Fibered modification) $\varphi$ is called a \emph{fibered map} if it restricts to an isomorphism of $\a^1$-fibered quasi-projective surfaces

$$\begin{xy}
\xymatrix{
V = X \backslash D \ar[r]^\sim_\varphi \ar[d]_{\bar{\pi}\mid_{V}} & V' = X' \backslash D' \ar[d]^{\bar{\pi}'\mid_{V'}} \\
\a^1 \ar[r]^\sim & \a^1.
}
\end{xy}$$

\noindent $\varphi$ is called \emph{fibered modification} if it is not an isomorphism.
\item[(2)] (Reversion) $\varphi$ is called \emph{reversion} if it admits a resolution of the form

$$\begin{xy}
\xymatrix{
 & (Z, \tilde{D} = C_n \triangleright \cdots \triangleright C_1 \triangleright H \triangleright C'_1 \triangleright \cdots \triangleright C'_{n'}) \ar[ld]_\sigma \ar[rd]^{\sigma'} & \\
(X, {}^tD) \ar@{-->}[rr]^\varphi &  & (X', D'),
}
\end{xy}$$

\noindent where $H$ is a zigzag with boundaries $C_0$ (left) and $C'_0$ (right) and where $\sigma: Z \to X$ and $\sigma': Z \to X'$ are smooth contractions of the sub-zigzags $H \triangleright C'_1 \triangleright \cdots \triangleright C'_{n'}$ and $C_n \triangleright \cdots \triangleright C_1 \triangleright H$ of $\tilde{D}$ onto $C_0$ and $C'_0$ respectively.
\end{itemize}
\end{df}

\begin{rem} We already introduced the notion of a reversion for standard pairs. Given a standard completion $(X, D)$ of $V$, we blow up $X$ in $C_0 \cap C_1$ and contract the proper transform of $C_0$. Letting $(X', D')$ be the resulting $1$-standard completion of $V$, we obtain a birational map $(X, D) \dasharrow (X', D')$. We will see in Proposition \ref{UniquenessReversions} below that these two notions of a reversion coincide after performing such elementary transformations on the boundary.
\end{rem}

It turns out that fibered modifications are just the liftings of appropriate triangular automorphisms of $\a^2$:

\begin{lemma}\label{JonquieresAutomorphisms}(\cite{BD}, Lemma 2.2.3) Let $\varphi: (X, D, \bar{\pi}) \dasharrow (X', D', \bar{\pi}')$ be a birational map between $1$-standard pairs and let $X \stackrel{\mu}{\leftarrow} \tilde{X} \stackrel{\eta}{\to} \f_1$ and $X' \stackrel{\mu'}{\leftarrow} \tilde{X}' \stackrel{\eta'}{\to} \f_1$ be as in Lemma \ref{FactorizationHirzebruch}. Then the following are equivalent:
\begin{itemize}
\item[(1)] $\varphi$ restricts to an isomorphism $(X \backslash D, \pi) \stackrel{\sim}{\to} (X' \backslash D', \pi')$.
\item[(2)] $(\mu')^{-1} \circ \varphi \circ \mu: \tilde{X} \dasharrow \tilde{X}'$ is the lift via $\eta$ and $\eta'$ of an isomorphism of affine $\a^1$-fibered surfaces

$$\begin{xy}
\xymatrix{
\a^2 = \f_1 \backslash (\eta(C_0) \cup \eta(C_1)) \ar[r]^\sim_\Psi \ar[d]_{\rho\mid_{\a^2}} & \a^2 = \f_1 \backslash (\eta'(C_0) \cup \eta'(C_1)) \ar[d]^{\rho\mid_{\a^2}} \\
\a^1 \ar[r]^\sim_{\psi} & \a^1,
}
\end{xy}$$

\noindent which maps isomorphically the base points of $\eta^{-1}$ onto those of $(\eta')^{-1}$. The map $\Psi$ is of the form $\Psi(x_0, y_0) = (ax_0 + P(y), by + c)$ with $a, b \in \c^*$, $c \in \c$ and $P(y) \in \c[y]$.
\end{itemize}
Moreover, $\varphi: (X, D) \dasharrow (X', D')$ is an isomorphism if and only if $\Psi$ is affine.
\end{lemma}

Indeed, the center $p$ of $\varphi$ gives the full control over the reversion:

\begin{prop}\label{UniquenessReversions} (Uniqueness of reversions, \cite{BD}, Prop. 2.3.7) For every $1$-standard pair $(X, D)$ and every point $p \in C_0 \backslash C_1$ there exist a $1$-standard pair $(X', D')$ and a reversion $\varphi: (X, D) \dasharrow (X', D')$, unique up to an isomorphism at the target, having $p$ as a unique proper base point. Moreover, if $\Gamma_D = [[0, -1, w_2, \dots, w_n]]$, then $\Gamma_{D'} = [[0, -1, w_n, \dots, w_2]]$.
\end{prop}

These two types of maps, the fibered modifications and the reversions, differ in the position of their center:

\begin{lemma}(\cite{BD}, Lemma 2.4.1) Let $\varphi: (X, D) \dasharrow (X', D')$ be a birational map between $1$-standard pairs.
\begin{itemize}
\item[(a)] If $\varphi$ is a fibered modification, it is centered at $p = C_0 \cap C_1$, and $C'_0$ is the only irreducible component of $D'$ contracted by $\varphi^{-1}$.
\item[(b)] If $\varphi$ is a reversion, it is centered at $p \in C_0 \backslash C_1$, and $\varphi^{-1}$ contracts the curves $C'_i$, $i \geq 1$, to $p$ as well as $C'_0$ to $p$ if and only if ${C'_i}^2 \leq -3$ holds for some $i \geq 2$.
\end{itemize}
\end{lemma}

If the type of the sub-zigzag $C_2 \triangleright \cdots \triangleright C_n$ is not a palindrome, then the composition of two reversions cannot be a reversion. Otherwise, we have the following

\begin{lemma}\label{CompositionOfReversions}(\cite{BD}, Lemma 2.3.8) For $i = 1, 2$ let $\varphi_i: (X, D) \dasharrow (X_i, D_i)$ be a reversion of $1$-standard pairs and assume that every component of $D$ has self-intersection $\geq -2$. If the proper base points of $\varphi_1$ and $\varphi_2$ are distinct (respectively equal), then the map $\varphi_2 \circ \varphi^{-1}_1$ is a reversion (respectively an isomorphism).
\end{lemma}

The key observation to control birational maps between $1$-standard pairs is to decompose any such map into fibered modifications and reversions:

\begin{thm}\label{FactorizationOfBirationalMaps}(\cite{BD}, Theorem 3.0.2) Let $\varphi: (X, D) \dasharrow (X', D')$ be a birational map between $1$-standard pairs restricting to an isomorphism $X \backslash D \stackrel{\sim}{\to} X' \backslash D'$. If $\varphi$ is not an isomorphism, then it can be decomposed into a finite sequence

$$\varphi = \varphi_n \circ \cdots \circ \varphi_1: (X, D) = (X_0, D_0) \stackrel{\varphi_1}{\to} (X_1, D_1) \stackrel{\varphi_2}{\to} \cdots \stackrel{\varphi_n}{\to} (X_n, D_n) = (X', D')$$

\noindent of fibered modifications and reversions between $1$-standard pairs $(X_i, D_i)$. Moreover, such a factorization of minimal length is unique, meaning, if 

$$\varphi = \varphi'_n \circ \cdots \circ \varphi'_1: (X, D) = (X'_0, D'_0) \stackrel{\varphi'_1}{\to} (X'_1, D'_1) \stackrel{\varphi'_2}{\to} \cdots \stackrel{\varphi'_n}{\to} (X'_n, D'_n) = (X', D')$$

\noindent is another factorization of minimal length, then there exist isomorphisms of $1$-standard pairs $\alpha_i: (X_i, D_i) \to (X'_i, D'_i)$, such that $\alpha_i \circ \varphi_i = \varphi'_i \circ \alpha_{i - 1}$ for $i = 2, \dots, n$.
\end{thm}

\subsection{Automorphisms of $\a^1$-fibrations and associated graphs}

Following \cite{DG2} and \cite{BD}, we introduce for an $\a^1$-fibered surface $V$ a (not necessarily finite) graph $\F_V$ which reflects the structure of the automorphism group of $V$.

\begin{df} To every normal affine surface $V$ we associate the oriented graph $\F_V$ as follows:
\begin{itemize}
\item[(1)] A vertex of $\F_V$ is an equivalence class of a $1$-standard pair $(X, D)$, such that $X \backslash D \cong V$, where two $1$-standard pairs $(X_1, D_1, \bar{\pi}_1)$ and $(X_2, D_2, \bar{\pi}_2)$ define the same vertex if and only if $(X_1 \backslash D_1, \pi_1) \cong (X_2 \backslash D_2, \pi_2)$.
\item[(2)] An arrow of $\F_V$ is an equivalence class of reversions. If $\varphi: (X, D) \to (X', D')$ is a reversion, then the class $[\varphi]$ of $\varphi$ is an arrow starting from $[(X, D)]$ and ending at $[(X', D')]$. Two reversions $\varphi_1: (X_1, D_1) \dasharrow (X'_1, D'_1)$ and $\varphi_2: (X_2, D_2) \dasharrow (X'_2, D'_2)$ define the same arrow if and only if there exist isomorphisms $\theta: (X_1, D_1) \to (X_2, D_2)$ and $\theta': (X'_1, D'_1) \to (X'_2, D'_2)$, such that $\varphi_2 \circ \theta = \theta' \circ \varphi_1$.
\end{itemize}
\end{df}

\begin{rem}\label{EquivalenceOfArrows} It follows from the definition that for a $1$-standard pair $(X, D)$ two reversions $\varphi_1: (X, D) \dasharrow (X_1, D_1)$ and $\varphi_2: (X, D) \dasharrow (X_2, D_2)$ centred at points $p_1$ and $p_2$ define the same arrow if and only if there exists an automorphism $\psi \in \Aut(X, D)$ such that $\psi(p_1) = p_2$.
\end{rem}

The structure of the graph $\F_V$ allows us to decide, whether the automorphism group $\Aut(V)$ of $V$ is generated by automorphisms of $\a^1$-fibrations. Here we say that $\varphi \in \Aut(V)$ is an \textit{automorphism of $\a^1$-fibrations} if there exists an $\a^1$-fibration $\pi: V \to \a^1$, such that $\varphi$ induces an isomorphism $\varphi: (V, \pi) \stackrel{\sim}{\to} (V, \pi)$.

\begin{prop}\label{SurfaceTree}(\cite{BD}, Prop. 4.0.7) Let $V$ be a normal affine surface with a non-empty graph $\F_V$. Then the following holds:
\begin{itemize}
\item[(1)] The graph $\F_V$ is connected.
\item[(2)] There is a natural bijection between the set of vertices of $\F_V$ and the isomorphism classes of $\a^1$-fibrations on $V$. 
\item[(3)] Let $(X, D)$ be a $1$-standard pair with $X \backslash D \cong V$ and let $D$ contain at least one curve with self-intersection $\leq -3$. Then there is a natural exact sequence

$$1 \to H \to \Aut(V) \to \Pi_1(\F_V) \to 1,$$

\noindent where $H$ is the (normal) subgroup of $\Aut(V)$ generated by all automorphisms of $\a^1$-fibrations and $\Pi_1(\F_V)$ is the fundamental group of the graph $\F_V$. In particular, the graph $\F_V$ is a tree if and only if $\Aut(V)$ is generated by automorphisms of $\a^1$-fibrations on $V$.
\end{itemize}
\end{prop}

One can obtain a better description of the group $\Aut(V)$ by introducing the notion of a graph of groups.

\begin{df} (cf. \cite{Se}) A \emph{graph of groups} is a pair $(\F, \G)$ such that $\F$ is an oriented graph and $\G$ consists of a family of vertex groups $\{ G_v \mid v \in V(\F) \}$ and a family of edge groups $\{ G_\sigma \mid \sigma \in E(\F) \}$ satisfying the following conditions:
\begin{itemize}
\item[(1)] For every edge it holds $G_\sigma = G_{\sigma^{-1}}$.
\item[(2)] For every edge $\sigma$ there are monomorphisms $\kappa_\sigma: G_\sigma \to G_{s(\sigma)}$ and $\lambda_\sigma: G_\sigma \to G_{t(\sigma)}$ such that $\lambda_\sigma = \kappa_{\sigma^{-1}}$. Here the index $t$ denotes the target variety of $\sigma$ and the index $s$ denotes the source variety of $\sigma$.
\end{itemize}
A \emph{path} in $(\F, \G)$ is a sequence $(g_0, \sigma_1, g_1, \dots, \sigma_r, g_r)$, where $g_i \in G_{v_i}$ and $v_0, \sigma_1, v_1, \dots, \sigma_r, v_r$ is a path in $\F$. The homotopy equivalence relation $\simeq$ is the equivalence relation generated by the elementary homotopy equivalence relations $(\sigma, \lambda_\sigma(h), \sigma^{-1}, (\kappa_\sigma(h))^{-1}) \simeq (1)$ with $1 \in G_{s(\sigma)}$ and $(g, \sigma, 1, \sigma^{-1}, g') \simeq (gg')$. If $v$ is a vertex of $\F$ then the homotopy classes of closed paths starting and ending in $v$ form a group under the concatenation $(\dots, g)(g', \dots) = (\dots, gg', \dots)$. We denote this group by $\pi_1(\F, \G, v)$ and call it the \emph{fundamental group} of $(\F, \G)$ in $v$.
\end{df}

We can equip $\F_V$ in a natural way with a structure of a graph of groups.

\begin{df}\label{DfGraph} Let $V$ be a normal quasi-projective surface and let $\F_V$ be its associated graph. Then $\F_V$ admits a structure of a graph of groups by the following choice:
\begin{itemize}
\item[(1)] For any vertex $v$ of $\F_V$, fix a $1$-standard pair $(X_v, D_v, \bar{\pi}_v)$ in the class $v$. The group $G_v$ is equal to $\Aut(X_v \backslash D_v, \pi_v)$.
\item[(2)] For any arrow $\sigma$ of $\F_V$, fix a reversion $r_\sigma: (X_\sigma, D_\sigma, \bar{\pi}_\sigma) \dasharrow (X'_\sigma, D'_\sigma, \bar{\pi'}_\sigma)$ in the class of $\sigma$ and also an isomorphism $\mu_\sigma: (X'_\sigma \backslash D'_\sigma, \pi'_\sigma) \to (X_{t(\sigma)} \backslash D_{t(\sigma)}, \pi_{t(\sigma)})$. Then the group $G_\sigma$ is equal to

$$\{ (\varphi, \varphi') \in \Aut(X_\sigma, D_\sigma) \times \Aut(X'_\sigma, D'_\sigma) \mid r_\sigma \circ \varphi = \varphi' \circ r_\sigma \}$$

and the monomorphisms $\kappa_\sigma: G_\sigma \to G_{s(\sigma)}$ and $\lambda_\sigma: G_\sigma \to G_{t(\sigma)}$ are given by $\kappa_\sigma((\varphi, \varphi')) = \mu_{\sigma^{-1}} \circ \varphi \circ \mu^{-1}_{\sigma^{-1}}$ and $\lambda_\sigma((\varphi, \varphi')) = \mu_\sigma \circ \varphi' \circ \mu^{-1}_\sigma$.
\end{itemize}
\end{df}

The first version of the following theorem was shown by Danilov and Gizatullin (\cite{DG2}, Theorem 5) and connects the structure of the graph of groups on $\F_V$ with the automorphism group of $V$:

\begin{thm}\label{AutFundamentalGroup} (\cite{DG2}, Theorem 5 and \cite{BD}, Theorem 4.0.11) Let $(X, D)$ be a $1$-standard pair such that $D$ admits at least one component with self-intersection $\leq -3$ and let $V := X \backslash D$. If $\F_V$ is equipped with a structure of a graph of groups as in Definition \ref{DfGraph} then the fundamental group of the graph of groups obtained is isomorphic to $\Aut(V)$.  
\end{thm}

\begin{rem}\label{ExTreeWithTwoVertices} Let $\F$ be the graph $v \ \bullet \stackrel{\sigma}{\longleftrightarrow} \bullet \ w$ and $\G = \{ (G_v, G_w), (G_\sigma) \}$. We can identify $G_\sigma$ via $\lambda_\sigma$ and $\kappa_\sigma$ respectively with subgroups of $G_v$ and $G_w$ respectively. It is a well-known result that $\pi_1(\F, \G, v)$ is isomorphic to the amalgamated product $G_v \star_{G_\sigma} G_w$.
\end{rem}

\subsection{The Matching Principle}

In the following we summarize the Matching Principle given in \cite{FKZ4}, section 3. We consider a standard completion $(X, D)$ of a smooth Gizatullin surface $V$ as well as the reversed completion $(X^\vee, D^\vee)$ with $D = C_0 \cup \cdots \cup C_n$ and $D^\vee = C^\vee_0 \cup \cdots \cup C^\vee_n$. We let $\Gamma_D = [[0, 0, w_2, \dots, w_n]]$ and we denote the corresponding extended divisors by $D_{\text{ext}}$ and $D^\vee_{\text{ext}}$, respectively. By inner elementary transformations we can move the pair of zeros to the right by several places. In this way we obtain, for every $t$, $2 \leq t \leq n + 1$, a new completion $(W, E)$ of $V$ with boundary divisor $[[w_2, \dots, w_{t - 1}, 0, 0, w_t, \dots, w_n]]$, \ie

$$E = C^\vee_n \cup \cdots \cup C^\vee_{n + 2 - t} \cup C_{t - 1} \cup C_t \cup \cdots \cup C_n,$$

\noindent if we identify $C_i \subseteq X$ and $C^\vee_j \subseteq X^\vee$ with their proper transforms in $W$. In particular, we can write $E$ as $E = D^{\geq t - 1} \cup D^{\vee \geq n + 2 - t}$ with new weights $C^2_{t - 1} = C^{\vee 2}_{n + 2 - t} = 0$. For brevity we let

$$t^\vee := n + 2 - t.$$

\noindent Moreover, there are natural isomorphisms

\begin{eqnarray*} 
W \backslash D^{\vee \geq t^\vee} &=& W \backslash(C^\vee_n \cup \cdots \cup C^\vee_{t^\vee}) \cong X \backslash (C_0 \cup \cdots \cup C_{t - 2}),\\
W \backslash D^{\geq t - 1} &=& W \backslash(C_{t - 1} \cup \cdots \cup C_n) \cong X^\vee \backslash (C^\vee_0 \cup \cdots \cup C^\vee_{t^\vee - 1}).\\
\end{eqnarray*}

\begin{df} (\cite{FKZ4}, Def. 3.3.3) The map

$$\psi := \Phi_{|C_{t - 1}|}: W \to \p^1$$

\noindent is called the \emph{correspondence fibration} for the pair $(C_t, C^\vee_{t^\vee})$.
\end{df}

There is a natural correspondence between feathers of $D_{\text{ext}}$ and those of $D^\vee_{\text{ext}}$. Before giving the key observation, we recall that for a given feather $F$ of $D_{\text{ext}}$, a boundary component $C_\mu$ is called \emph{mother component of} $F$ if the feather $F$ is created by a blowup on $C_\mu$ during the blowup process $X \to Q = \p^1 \times \p^1$ (see \cite{FKZ3}, 2.3).

\begin{prop}\label{MatchingFeathers} (cf. \cite{FKZ4}, Lemma 3.3.4, Cor. 3.3.5, Lemma 3.3.6) Let $F_{i, \rho}$ be a feather of $D_{\text{ext}}$ attached to the component $C_i$. Then there exists a unique feather $F^\vee_{j, \sigma}$ of $D^\vee_{\text{ext}}$ which intersects $F_{i, \rho}$ in $V$ and which is attached to a component $C^\vee_j$, such that $i + j \geq n + 2$. Moreover, $F_{i, \rho}$ and $F^\vee_{j, \sigma}$ intersect transversally and in a single point. If $C_\tau$ is the mother component of $F_{i, \rho}$, then $C^\vee_{\tau^\vee}$ is the mother component of $F^\vee_{j, \sigma}$. 
\end{prop}

\begin{df} Feathers $F_{i, \rho}$ and $F^\vee_{j, \sigma}$, which satisfy the conditions of Proposition \ref{MatchingFeathers}, are called \emph{matching feathers}.
\end{df}

The condition $i + j \geq n + 2$ is essential. Indeed, every feather $F_{t - 1, \rho}$ is a section of $\psi$ and therefore it meets every fiber of $\psi$. Since it cannot intersect $D^{\geq t}_{\text{ext}}$, it meets every feather $F^\vee_{t^\vee, \sigma}$ with $(F^\vee_{t^\vee, \sigma})^2 = -1$ on $V$.

\bigskip
\noindent \textbf{Configuration spaces and the configuration invariant.} Let $V$ be a smooth Gizatullin surface with a standard completion $(X, D)$. The sequence of weights $[[w_2, \dots, w_n]]$ (up to reversion) of the boundary divisor $D$ is a discrete invariant of the abstract isomorphism type of $V$ (\cite{FKZ1C}, Cor. 3.33')). In the following we recall a stronger continuous invariant of $V$, the \textit{configuration invariant} (see \cite{FKZ4}, 3.2).\\
For a natural number $s \geq 1$ we denote the configuration space of all $s$-points subsets $\{ \lambda_1, \dots, \lambda_s \} \subseteq \a^1$ by $\mathcal{M}^+_s$. We can identify $\mathcal{M}^+_s$ in a natural way with the Zariski open subset of $\a^s$:

$$\mathcal{M}^+_s \cong \a^s \backslash \{ \text{discr}(P) = 0 \}, \quad \text{where} \quad P = \prod_{j = 1}^s (X - \lambda_j),$$

\noindent see \cite{FKZ4}, 3.1.1. The group $\Aut(\a^1)$ acts on $\mathcal{M}^+_s$ in a natural way. We let

$$\M^+_s := \mathcal{M}^+_s/\Aut(\a^1).$$

\noindent Thus, $\M^+_s$ is an $(s - 2)$-dimensional affine variety. Now, let $\mathcal{M}^*_s$ be the configuration space of all $s$-points subsets $\{ \lambda_1, \dots, \lambda_s \} \subseteq \c^* = \a^1 \backslash \{ 0 \}$. Similarly, the group $\c^*$ acts on $\mathcal{M}^*_s$ and we let 

$$\M^*_s := \mathcal{M}^*_s/\c^*.$$

Before introducing the configuration invariant we have to distinguish two types of boundary components.

\begin{df}\label{StarComponent}
\begin{itemize}
\item[(1)] For a natural number $i \in \{ 2, \dots, n \}$ $s_i$ shall denote the number of feathers of $D_{\text{ext}}$ whose mother component is $C_i$.
\item[(2)] The component $C_i$ is called a \emph{$*$-component} or \emph{inner component} if 
\begin{itemize}
\item[(i)] $D^{\geq i + 1}_{\text{ext}}$ is not contractible and
\item[(ii)] $D^{\geq i + 1}_{\text{ext}} - F_{j, k}$ is not contractible for every feather $F_{j, k}$ of $D^{\geq i + 1}_{\text{ext}}$ with mother component $C_\tau$, where $\tau < i$.
\end{itemize}
Otherwise $C_i$ is called a \emph{$+$-component} or \emph{outer component}.
\end{itemize}
\end{df}

For example, $C_2$ and $C_n$ are always $+$-components. In the following we let $\tau_i = *$ in the first case and $\tau_i = +$ in the second one.\\
It is easily seen that in the blowup process $\tilde{X} \to \p^1 \times \p^1$ ($\tilde{X}$ is a standard completion of the minimal resolution of singularities $V'$ of $V$) every $*$-component $C_i$, $3 \leq i \leq n - 1$, appears as a result of an inner blowup of the previous zigzag, while an outer blowup of a zigzag creates a $+$-component.\\
The following lemma states that reversions do not change the type of a component.

\begin{lemma}\label{ComponentTypeAfterReversion} (\cite{FKZ4}, Lemma 3.3.10) $C_t$ is a $*$-component if and only if $C^\vee_{t^\vee}$ is a $*$-component.
\end{lemma}

Now we are able to construct the so-called \textit{configuration invariant of} $V$. In the following, for a $+$-component $C_i$ we construct a family of points $p_{i, j}$, $1 \leq j \leq s_i$, on $C_i \backslash C_{i - 1} \cong \a^1$. For every feather $F_{i, j}$ with self-intersection $-1$ we let $p_{i, j}$ be its intersection point with $C_i$. Moreover, if there exists a feather $F_{k, j}$ with mother component $C_i$ and $k > i$, then we also add the intersection point $c_{i + 1} := C_i \cap C_{i + 1}$ to our collection (note, that such a feather is unique, if it exists). Thus, the collection of points

$$p_{i, j} \in C_i, \quad 1 \leq j \leq s_i$$

\noindent is just the collection of locations on $C_i$ in which the feathers with mother component $C_i$ are born by a blowup. These points are called \textit{base points} of the associated feathers. The collection $(p_{i, j})_{1 \leq j \leq s_i}$ defines a point $Q_i$ in $\M^+_{s_i}$.

Let now $C_i$ be a $*$-component. Then we consider $Q_i$ as a collection of points on $C_i \backslash (C_{i - 1} \cup C_{i + 1})$. Note that the intersection point $c_{i + 1}$ of $C_i$ and $C_{i + 1}$ cannot belong to this collection due to Definition \ref{StarComponent} (2) (ii). Identifying $C_i \backslash (C_{i - 1} \cup C_{i + 1})$ with $\c^*$ in a way that $c_{i + 1}$ corresponds to $0$ and $c_i$ to $\infty$ we obtain a point $Q_i$ in the configuration space $\M^*_{s_i}$.\\
Thus, in total, we obtain a point 

$$Q(X, D) := (Q_2, \dots, Q_n) \in \M = \M^{\tau_2}_{s_2} \times \cdots \times \M^{\tau_n}_{s_n},$$

\noindent where $\tau_i \in \{ +, * \}$ represents the type of the corresponding component $C_i$. This point is called the \textit{configuration invariant of} $(X, D)$.
 
Performing elementary transformations in $(X, D)$ with centers in $C_0$ does neither change $\Phi_0$ nor the extended divisor (except for the weight $C_1^2$). Thus, it leaves the $s_i$ and $Q(X, D)$ invariant. Therefore, we can define the configuration invariant for every $m$-standard completion of $V$.\\

\begin{prop}\label{MatchingPrinciple} (Matching Principle, \cite{FKZ4}, Prop. 3.3.1) Let $V = X \backslash D$ be a smooth Gizatullin surface completed by a standard zigzag $D$. Consider the reversed completion $(X^\vee, D^\vee)$ with boundary zigzag $D^\vee = C_0^\vee \cup \cdots \cup C_n^\vee$, the associated numbers $s'_2, \dots, s'_n$ and the types $\tau'_2, \dots, \tau'_n$. Then $s_i = s'_{i^\vee}$ and $\tau_i = \tau'_{i^\vee}$ for all $i = 2, \dots, n$. Moreover, the associated points $Q(X, D)$ and $Q(X^\vee, D^\vee)$ in $\M$ coincide under the natural identification

$$\M = \M^{\tau_2}_{s_2} \times \cdots \times \M^{\tau_n}_{s_n} \cong \M^{\tau'_n}_{s'_n} \times \cdots \times \M^{\tau'_2}_{s'_2}.$$

\end{prop}

\begin{df} Let $\Gamma$ and $\Gamma'$ be weighted graphs. A \emph{reconstruction} $\gamma$ of $\Gamma$ into $\Gamma'$ is a finite sequence

$$\gamma: \Gamma = \Gamma_0 \stackrel{\gamma_1}{\dasharrow} \Gamma_1 \stackrel{\gamma_2}{\dasharrow} \cdots \stackrel{\gamma_n}{\dasharrow} \Gamma_n = \Gamma',$$

\noindent where each arrow $\gamma_i$ is either a blowup or a blowdown. The graph $\Gamma'$ is called \emph{end graph} of $\gamma$. The inverse sequence $\gamma^{-1} := (\gamma_n^{-1}, \dots, \gamma_1^{-1})$ yields a reconstruction of $\Gamma'$ with end graph $\Gamma$.\\
The reconstruction $\gamma$ is said to be \emph{symmetric} if it is of the form $(\gamma, \gamma^{-1})$.
\end{df}

\begin{df} Two standard completions $(X, D)$ and $(X', D')$ of a Gizatullin surface $V$ are \emph{evenly linked} if there is a symmetric reconstruction of $(X, D)$ into $(X', D')$. Otherwise they are called \emph{oddly linked}.
\end{df}

Indeed, a standard completion of a Gizatullin surface is evenly linked to any other standard completion or to its inverse:

\begin{lemma}\label{EvenlyLinked} (\cite{FKZ4}, Lemma 2.2.2) Let $(X, D)$ and $(X', D')$ be two standard completions of a Gizatullin surface $V \not\cong \a^1 \times \c^*$. After replacing, if necessary, $(X, D)$ by its reversion $(X^\vee, D^\vee)$, $(X, D)$ and $(X', D')$ are evenly linked.
\end{lemma}

The next theorem shows that the configuration invariant $Q(V)$ of a smooth Gizatullin surface $V$ is indeed an invariant of the abstract isomorphism type of $V$:

\begin{thm}\label{InvarianceConfigInvariant} (\cite{FKZ4}, Theorem 3.4.1) Given two $m$-standard completions $(X, D), (X', D')$ of a smooth Gizatullin surface $V$, for the configuration invariants $s_i, s'_i$ and $Q(X, D) \in \M$, $Q(X', D') \in \M'$ the following holds:
\begin{itemize}
\item[(1)] If $(X, D)$ and $(X', D')$ are evenly linked, then $s_i = s'_i$ for $i = 2, \dots, n$ and the points $Q(X, D)$ and $Q(X', D')$ of $\M = \M'$ coincide.
\item[(2)] If $(X, D)$ and $(X', D')$ are oddly linked, then $s_i = s'_{i^{\vee}}$ for $i = 2, \dots, n$ and the points $Q(X, D) \in \M$ and $Q(X', D') \in \M'$ of $\M$ and $\M'$ coincide under the natural identification

$$\M = \M^{\tau_2}_{s_2} \times \cdots \times \M^{\tau_n}_{s_n} \cong \M^{\tau'_n}_{s'_n} \times \cdots \times \M^{\tau'_2}_{s'_2} = \M'.$$

\end{itemize}
\end{thm}

\begin{df} Given a configuration space $\M = \M^{\tau_2}_{s_2} \times \cdots \times \M^{\tau_n}_{s_n}$ we consider the reversed product

$$\M^\vee := \M^{\tau_n}_{s_n} \times \cdots \times \M^{\tau_2}_{s_2}.$$

\noindent The \emph{symmetric configuration invariant} of a completion $(X, D)$ of a smooth Gizatullin surface $V$ is the unordered pair

$$\tilde{Q}(X, D) := \{ Q(X, D), Q(X^\vee, D^\vee) \}, \quad \text{where} \quad Q(X, D) \in \M \quad \text{and} \quad Q(X^\vee, D^\vee) \in \M^\vee.$$

\end{df}

Now, the following is obvious:

\begin{cor} (\cite{FKZ4}, Cor. 3.4.3) The pair $\tilde{Q}(V) := \tilde{Q}(X, D)$ as well as the sequence $(s_2, \dots, s_n)$ (up to reversion) are invariants of the isomorphism type of $V$.
\end{cor}

\subsection{Coordinates on smooth Gizatullin surfaces}\label{CoordinatesGeneralCase}

For our purpose we need explicit descriptions of smooth Gizatullin surfaces via affine coordinates on appropriate open affine charts. In \cite{FKZ4}, \S 4, such coordinates are constructed for the case of smooth Gizatullin surfaces which admit a \emph{presentation}, \ie where all boundary components $C_i$, $i \geq 2$ are $+$-components. We generalize this description for the case that a smooth Gizatullin surface may have inner components. Our interest will be concentrated on surfaces admitting a standard completion $(X, D)$ such that $C_3, \dots, C_{n - 1}$ are inner components and such that no feathers are attached to $C_2$ and $C_n$.
 
Let study the process which creates the $*$-components $C_3, \dots, C_{n - 1}$ by successive inner blowups. A standard completion $(X, D)$ of a smooth Gizatullin surface $V$ can be realized as a sequence of blowups of the quadric $Q = \p^1 \times \p^1$, such that all blowup centers are contained in $C_2 \backslash C_1$ and its infinitely near neighbourhood. In contrast to the case where all $C_i$ are outer components, there is no prescribed ordering for creating inner components. However, the algorithm below describes local coordinate charts on $X$, and the ordering does not play any role in our final results.

We consider on $X_0 := Q = \p^1 \times \p^1$ the affine chart $Q \backslash (C_0 \cup C_1) \cong \a^2$ with affine coordinates $(x_0, y_0)$, such that $C_0 = \{ y_0 = \infty \}$, $C_1 = \{ x_0 = \infty \}$ and $C_2 = \{ y_0 = 0 \}$ holds and decompose the map $X \to \p^1 \times \p^1$ into blowups

$$X = X_N \stackrel{\pi_N}{\to} X_{N - 1} \stackrel{\pi_{N - 1}}{\to} \cdots \stackrel{\pi_1}{\to} X_0 = Q,$$

\noindent where every $\pi_i$ creates either a new boundary component or a family of feathers attached to the same component. Now we proceed as follows:\\

(1) Without loss of generality we may assume that $C_2 \cap C_3$ has coordinates $(x_0, y_0) = (0, 0)$. We let $(s_0, t_0) := (x_0, y_0)$ and introduce affine coordinates inductively either via $(s_i, t_i) = (s_{i + 1}, t_{i + 1})$ (this corresponds to the case when the blowup is performed at infinity) or via $(s_i, t_i) = (s_{i + 1}t_{i + 1}, t_{i + 1})$ or via $(s_i, t_i) = (s_{i + 1}, s_{i + 1}t_{i + 1})$. Note, that these local coordinates have the property that both axes are given by certain boundary components. After $n - 2$ steps we arrive at $C_k \backslash C_{k - 1} = \{ t_{n - 2} = 0 \}$ and $C_{k + 1} \backslash C_{k + 2} = \{ s_{n - 2} = 0 \}$ for some $k$ (cf. the following figure for the case of the second transformation; in this figure the proper transform of $C_{j + 1}$ becomes $C_{j + 2}$, since the exceptional curve of the blowup precedes the proper transform of $C_{j + 1}$). We let $(u_k, v_k) := (s_{n - 2}, t_{n - 2})$.

\bigskip
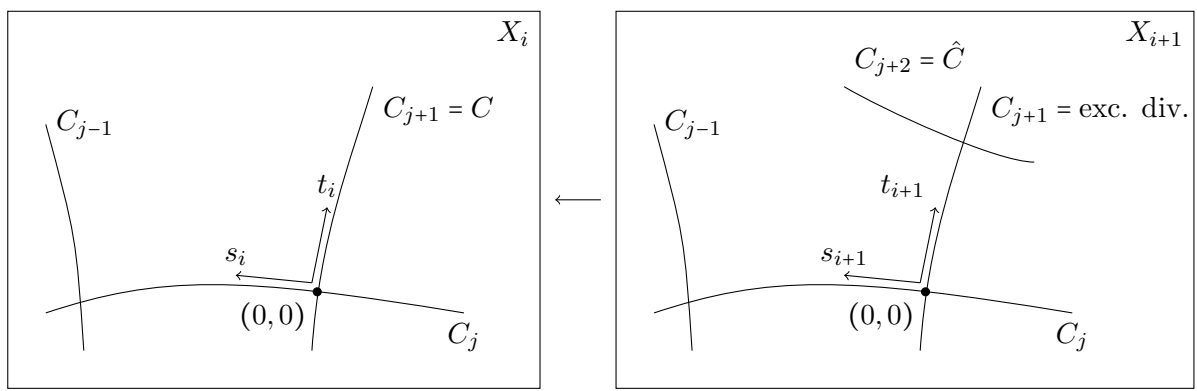
\begin{tikzpicture}
	\draw (0, 0) rectangle (7, 5) node[anchor=north east]{$X_i$};
	\draw (0.5, 3.5) node[anchor=west] {$C_{j - 1}$} .. controls (0.9, 2) and (0.9, 2) .. (1, 0.5); 			
	\draw (0.5, 1) .. controls (2, 1.5) and (3, 1.5) .. (6, 1) node[anchor=north] {$C_j$};					
	\draw (4, 0.5) .. controls (4.1, 2) and (4.5, 3) .. (4.8, 4) node[anchor=north west] {$C_{j + 1} = C$};					
	\filldraw [black] (4.07, 1.28) node[anchor=north east] {$(0, 0)$} circle (1.5pt);
	\draw [->] (4, 1.4) -- (3, 1.5) node[anchor=south] {$s_i$};						
	\draw [->] (4, 1.4) -- (4.2, 2.4) node[anchor=south] {$t_i$};						

	\draw [->] (7.8, 2.5) -- (7.2, 2.5);																		

	\draw (8, 0) rectangle (15.6, 5) node[anchor=north east]{$X_{i + 1}$};
	\draw (8.5, 3.5) node[anchor=west] {$C_{j - 1}$} .. controls (8.9, 2) and (8.9, 2) .. (9, 0.5); 			
	\draw (8.5, 1) .. controls (10, 1.5) and (11, 1.5) .. (14, 1) node[anchor=north] {$C_j$};					
	\draw (12, 0.5) .. controls (12.1, 2) and (12.5, 3) .. (12.8, 4) node[anchor=north west] {$C_{j + 1} = \text{exc. div.}$};					
	\draw (13.5, 3) .. controls (13, 3) and (11.5, 3.7) .. (11, 4) node[anchor=south west] {$C_{j + 2} = \hat{C}$};					
	\filldraw [black] (12.07, 1.28) node[anchor=north east] {$(0, 0)$} circle (1.5pt);
	\draw [->] (12, 1.4) -- (11, 1.5) node[anchor=south] {$s_{i + 1}$};						
	\draw [->] (12, 1.4) -- (12.2, 2.4) node[anchor=south east] {$t_{i + 1}$};						
\end{tikzpicture}

\bigskip
\noindent Now, in order to preserve the description of the correspondence fibration and, in particular, the one of the intersection points $F_{i + 1} \cap F^\vee_{i + 1}$ (cf. \cite{FKZ4}, 5.1, 5.2), we proceed as follows.

(2) Let the first $n - 2$ blowups $X_{n - 2} \to \cdots \to X_0 = Q$ create all inner components. It is easy to check that the coordinate systems $(u_i, v_i)$, $2 \leq i \leq n - 1$, from step (1) satisfy certain relations

$$(u_j, v_j) = T_{ij}(u_i, v_i) := (u_i^{k_{ij}}v_i^{l_{ij}}, u_i^{p_{ij}}v_i^{q_{ij}}) \quad \text{with} \quad \left|\begin{array}{cc} k_{ij} & p_{ij} \\ l_{ij} & q_{ij} \end{array}\right| = 1 \quad \text{and} \quad q_{ij} > 0, \ l_{ij} < 0 \quad \forall \ j > i.$$

(3) Given a finite subset $M \subseteq \c^*$, we let $P_M(x) := \prod_{a \in M} (x - a) \in \c[x]$. For $i = 3, \dots, n - 1$ we consider the base points $(a_{i, 1}, 0), \dots, (a_{i, r_i}, 0) \in C_i \backslash (C_{i - 1} \cup C_{i + 1})$ (in coordinates $(u_i, v_i)$) of the feathers attached to $C_i$. Letting $X_s$ be the underlying intermediate surface, we introduce new affine coordinates $(x_i, y_i)$ after a blowup $X_{s + 1} \to X_s$ in $M_i$ via

$$(x_i, y_i) := T^{\text{Bl.up}}_k(u_i, v_i) := \left( u_i, \frac{v_i}{P_{M_i}(u_i)} \right), \quad M_i := \{ a_{i, 1}, \dots, a_{i, r_i} \}.$$

\noindent Further, we let $F^\vee_{i, j}$, $j \in \{ 1, \dots, r_i \}$, be the proper transform of the closure of the affine line $\{ u_i = a_{i, j} \}$ under the blowup of $M_i$.

Now, using the relations $T_{ij}$ in (2) we replace the coordinates $(u_j, v_j)$ by 

$$(T_{ij} \circ T^{\text{Bl.up}}_i \circ T_{ij}^{-1})(u_j, v_j) \quad \text{for all} \quad j \in \{ i + 1 \dots, n \}.$$

\vspace{20pt}
\noindent Indeed, the following Proposition shows that the projection $(x_i, y_i) \mapsto x_i$ gives the correspondence fibration for the pair $(C^\vee_{j^\vee}, C_j)$:

\begin{prop}\label{CorrespondenceFibration} Assume that there are given affine coordinates $(x_i, y_i)$ on $X_N$ as above. Then, in appropriate coordinates on $\a^1 = \p^1 \backslash \{ \infty \}$, the map $\Phi_{|C^\vee_{j^\vee - 1}|}: W \to \p^1$, being the correspondence fibration for the pair $(C^\vee_{j^\vee}, C_j)$, is given by $x_j$. In particular, the pair $(F_{i, j}, F^\vee_{i, j})$ is a pair of matching feathers.
\end{prop}

\begin{bew} First of all, we note that if $(u_j, v_j)$ are local coordinates on $X$, then $(u'_j, v'_j) := (T_{ij} \circ T^{\text{Bl.up}}_i \circ T_{ij}^{-1})(u_j, v_j)$ also give well-defined local coordinates on $X$ (in the sense that for every given $(u'_j, v'_j) \in \c^2$ there is a unique point $x \in X$ with such coordinates). Hence, every coordinate system occuring in the algorithm is well-defined.\\

We show the claim by induction on $j \in \{ 2, \dots, n - 1 \}$. For $j = 2$, the caim is obvious. So let us assume that $(x_j, y_j) \mapsto x_j$, $j \geq 2$, gives the correspondence fibration for $(C^\vee_{j^\vee}, C_j)$. To describe the correspondence fibration for the pair $(C^\vee_{{j + 1}^\vee}, C_{j + 1})$ it is sufficient to construct the intermediate surface $X_M$, $M \leq N$, where only feathers attached to $C_\mu$, $\mu \leq j + 1$, are already created. First, let us consider two further intermediate surfaces $X_{M'}$ and $X_{M''}$, $M'' \leq M' \leq M$, such that

$\bullet$ On $X_{M'}$ only feathers which are attached to $C_\mu$, $\mu \leq j$, are already created. 

$\bullet$ On $X_{M''}$ only feathers which are attached to $C_\mu$, $\mu \leq j - 1$, are already created. 

\noindent We partially reverse $(X_{M''}, D_{M''})$ until the part ${D^\vee}^{\geq j^\vee}$ is constructed. On the resulting surface $(W'', E'')$, the curves $C^\vee_{j^\vee}$ and $C_{j - 1}$ are $0$-curves. Let further $W'$ be the surface corresponding to the correspondence fibration for $(C_{j + 1}, C^\vee_{(j + 1)^\vee})$, constructed from $X_{M'}$. Since $C^\vee_{j^\vee - 1}$ is the general fiber of the correspondence fibration for the pair $(C^\vee_{(j + 1)^\vee}, C_{j + 1})$, it is sufficient to show that the map $(x_{j + 1}, y_{j + 1}) \to x_{j + 1}$ extends to a map $W' \to \p^1$, such that its general fiber is isomorphic to $\p^1$ and such that $C^\vee_{{j + 1}^\vee - 1} = \overline{\{ y_j = 0 \}}$. Indeed, if two $\p^1$-fibrations coincide in a general fiber $C$, then they coincide everywhere, since both fibrations are given by the morphism associated to the linear system $|C|$.

Now, on $(W'', E'')$ we have three coordinate systems, $(x_{j - 1}, y_{j - 1})$, $(u_j, v_j)$ and $(u_{j + 1}, v_{j + 1})$, such that $C_{j - 1} = \overline{\{ y_{j - 1} = 0 \}}$, $C_j = \overline{\{ x_{j - 1} = 0 \}} = \overline{\{ y_j = 0 \}}$, $C_{j + 1} = \overline{\{ x_j = 0 \}} = \overline{\{ v_{j + 1} = 0 \}}$ and $C_{j + 2} = \overline{\{ u_{j + 1} = 0 \}}$.

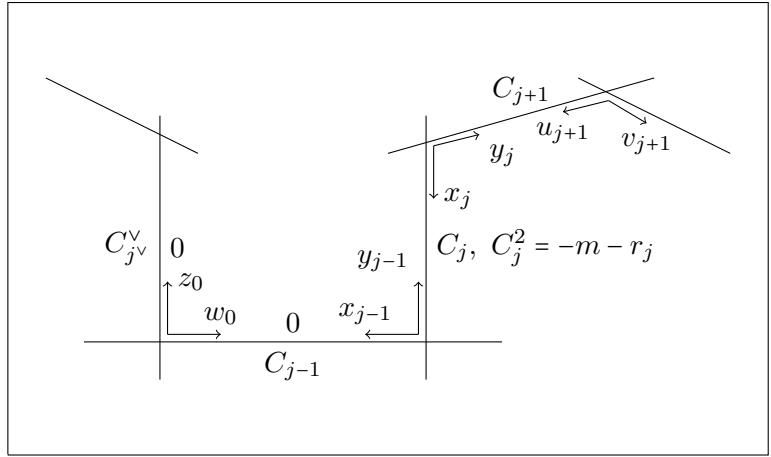
\vspace{20pt}
\bigskip
\begin{tikzpicture}
	\draw (3, 0) rectangle (13, 6);
	\begin{scope}
		\draw[color=white] (0, 0) -- (3, 0);
		
		\draw (3.5, 5) -- (5.5, 4);				

		\draw (5, 4.5) -- node[anchor=east] {$C^\vee_{j^\vee}$} 								node[anchor=west] {$0$} (5, 1);	
		\draw (4, 1.5) -- node[anchor=north] {$C_{j - 1}$} node[anchor=south] 					{$0$} (9.5, 1.5);	
		\draw (8.5, 4.5) -- node[anchor=west] {$C_j, \ C^2_j = -m - r_j$} (8.5, 1);					
		\draw (8, 4) -- node[anchor=south] {$C_{j + 1}$} (11.5, 5);				
		\draw (10.5, 5) -- (12.5, 4);				

		\draw [->] (5.1, 1.6) -- (5.8, 1.6) node[anchor=south] {$w_0$};	
		\draw [->] (5.1, 1.6) -- (5.1, 2.3) node[anchor=west] {$z_0$};	

		\draw [->] (8.4, 1.6) -- (7.7, 1.6) node[anchor=south] {$x_{j - 1}$};	
		\draw [->] (8.4, 1.6) -- (8.4, 2.3) node[anchor=south east] {$y_{j - 1}$};	

		\draw [->] (8.6, 4.1) -- (8.6, 3.4) node[anchor=west] {$x_j$};	
		\draw [->] (8.6, 4.1) -- (9.2, 4.25) node[anchor=north west] {$y_j$};	

		\draw [->] (10.9, 4.7) -- (10.3, 4.55) node[anchor=north] {$u_{j + 1}$};	
		\draw [->] (10.9, 4.7) -- (11.4, 4.4) node[anchor=north] {$v_{j + 1}$};	
		
	\end{scope}	
\end{tikzpicture}

\vspace{20pt}
\noindent By step (2) and (3) of the algorithm, which replace $(u_j, v_j)$ by $(x_j, y_j)$, it follows that

$$(u_j, v_j) = \left( \frac{1}{y_{j - 1}}, x_{j - 1}y_{j - 1}^m \right),$$

\noindent where $-m$ is the self-intersection number of the proper transform of $C_j$ on the surface, which is obtained by blowing down all feathers on $C_j$ (in other words, we have $C^2_j = -m - r_j$). Let now $\{ (a_1, 0), \dots, (a_{r_j}, 0) \} \subseteq \c^* \cong C_j \backslash (C_{j - 1} \cup C_{j + 1})$ be the $(u_j, v_j)$-coordinates of the base points of the feathers $F_{j, k}$, $k \in \{ 1, \dots, r_j \}$. Creating these feathers (which results in the surface $W'$), step (3) gives us

$$(x_j, y_j) = \left( u_j, \frac{v_j}{\prod_1^{r_j} (u_j - a_{j, i})} \right), \quad \text{or equivalently} \quad (u_j, v_j) = \left( x_j, y_j\prod_1^{r_j} (x_j - a_{j, i}) \right),$$

\noindent hence

\begin{equation}\label{CoordinatesForLines}
(x_{j - 1}, y_{j - 1}) = \left( x_j^my_j\prod_1^{r_j} (x_j - a_{j, i}), \frac{1}{x_j} \right).
\end{equation}

We consider an affine line $L_a$ with equation $y_j = a$ for some $a \in \c^*$. We show that its closure $\bar{L}_a$ intersects $C^\vee_{j^\vee - 1}$, and that these projective lines do not intersect for different values of $a$. This gives us, that the $\a^1$-fibration $(x_j, y_j) \mapsto y_j$ extends to a $\p^1$-fibration on $W'$, such that $C^\vee_{j^\vee - 1}$ is a section. Let us parametrize $L_a$ via $(x_j, y_j) = (t, a)$, $t \in \c$. We introduce local coordinates $(w_0, z_0)$, centred at $P := C^\vee_{j^\vee} \cap C_{j - 1}$, via $(w_0, z_0) = (1/x_{j - 1}, y_{j - 1})$. Then (\ref{CoordinatesForLines}) gives us, that the equation of $L_a$ is given by

$$L_a = \left\{ (w_0, z_0) = \left( \frac{1}{at^m\prod_1^{r_j} (t - a_{j, i})}, \frac{1}{t} \right) \mid t \in \c \right\}.$$

\noindent Hence, letting $s := 1/t$, we obtain that

$$\bar{L}_a \backslash C_{j + 1} = \left\{ (w_0, z_0) = \left( \frac{1}{a} \cdot \frac{s^{m + r_j}}{\prod_1^{r_j} (1 - a_{j, i}s)}, s \right) \mid s \in \c \right\}.$$

\noindent Now, to obtain the coorespondence fibration for $(C^\vee_{j^\vee}, C_j)$, we have to perform $m + r_j$ elementary transformations centred at $P$. This leads to new local coordinates $(w, z) = (w_0/z_0^{m + r_j}, z_0)$ (with $C^\vee_{j^\vee} = \overline{\{ w = 0 \}}$ and $C^\vee_{j^\vee - 1} = \overline{\{ z = 0 \}}$). Thus, on the resulting surface we have

$$\bar{L}_a \backslash C_{j + 1} = \left\{ (w, z) = \left( \frac{1}{a} \cdot \frac{1}{\prod_1^{r_j} (1 - a_{j, i}s)}, s \right) \mid s \in \c \right\},$$

\noindent Hence $L_a$ intersects $C^\vee_{j^\vee - 1}$ in the point $(w, z) = (1/a, 0)$. It follows that the projective lines $\bar{L}_a$ intersect $C^\vee_{j^\vee - 1}$ and that they do not intersect for different values of $a$.

Now, on the surface $W'$, we have $u_{j + 1} = 1/y_j$, hence lines with constant $y_j$-coordinate have also constant $u_{j + 1}$-coordinate. Thus the claim follows.

Finally, creating the feathers $F_{j + 1, i}$, $i = 1, \dots, r_{j + 1}$, does neither change the general fibers for the correspondence fibration for $(C_{j + 1}, C^\vee_{(j + 1)^\vee})$, nor the coordinate $u_{j + 1}$. The same holds for the elementary transformation for passing to the correspondence fibration for $(C^\vee_{(j + 1)^\vee}, C_{j + 1})$. Now, since $C_j$ (and hence $C^\vee_{j^\vee}$) has the equation $y_j = 0$, it is a fiber of $(x_{j + 1}, y_{j + 1}) \mapsto x_{j + 1}$ (at infinity). Thus, $(x_{j + 1}, y_{j + 1}) \mapsto x_{j + 1}$ gives the correspondence fibration for the pair of curves $(C^\vee_{(j + 1)^\vee}, C_{j + 1})$.
\end{bew}

\section{Actions of the automorphism groups of Gizatullin surfaces}

In the following we work with standard pairs as well as with $1$-standard pairs. Let us introduce some notations concerning $1$-standard pairs:

\bigskip
\noindent \textbf{Notation:} According to Lemma 1.0.7 in \cite{BD}, any smooth $1$-standard pair $(X, D)$ may be obtained by some blowups of points on a fiber of $\f_1$. An embedding of $\f_1$ into $\p^2 \times \p^1$ is given by

$$\f_1 = \text{Bl}_{(1: 0: 0)}(\p^2) = \{ ((x: y: z), (s: t)) \in \p^2 \times \p^1 \mid yt - zs = 0 \}.$$

\noindent We denote by $\tau : \f_1 \to \p^2$ the projection on $\p^2$ and by $C_0$ and $C_2$ the lines $\{ z = 0 \}$ and $\{ y = 0 \}$ respectively. We also denote by $C_0$ and $C_2$ their proper transforms on $\f_1$, by $C_1$ the exceptional curve $\tau^{-1}(1: 0: 0) = \{(1: 0: 0)\} \times \p^1$ and by $L_0$ the affine line $C_2 \backslash C_1 \subseteq \f_1$ as well as its image $C_2 \backslash \{ (1: 0: 0) \}$ in $\p^2$. Moreover, we have isomorphisms

$$\a^2 \stackrel{\cong}{\to} \f_1 \backslash (C_0 \cup C_1),\quad (w_0, z_0) \mapsto ((w_0: z_0: 1), (z_0: 1))$$

\noindent as well as 

$$\a^2 \stackrel{\cong}{\to} \p^2 \backslash C_0,\quad (w_0, z_0) \mapsto (w_0: z_0: 1).$$

\noindent In these coordinates the affine line $L_0$ is given by $z_0 = 0$. In the following we will denote these coordinates on $\f_1 \backslash (C_0 \cup C_1)$ by $(w_0, z_0)$. The $\p^1$-fibration on $\f_1$ is given via the second projection:

$$\rho: \f_1 \to \p^1,\quad ((x: y: z), (s: t)) \mapsto (s: t),$$

\noindent The restriction of $\rho$ to $\a^2 = \f_1 \backslash (C_0 \cup C_1)$ yields an $\a^1$-fibration $\pi: \a^2 \to \a^1$, which is simply the projection $(w_0, z_0) \mapsto z_0$ onto the second factor.\\

We denote by $\Aff$ the group of automorphisms of $\a^2$, which extend to automorphisms of $\p^2$ and by $\Jon$ the group of triangular (de Jonquieres) automorphisms (automorphisms of $(\a^2, \pi)$). In other words, we have

\begin{eqnarray*}
\Aff &=& \{ (w_0, z_0) \mapsto (a_{11}w_0 + a_{12}z_0 + b_1, a_{21}w_0 + a_{22}z_0 + b_2) \mid a_{11}a_{22} - a_{12}a_{21} \neq 0 \} \\
\Jon &=& \{ (w_0, z_0) \mapsto (aw_0 + P(z_0), bz_0 + c) \mid a, b \in \c^*, c \in \c, P(z_0) \in \c[z_0] \}.
\end{eqnarray*}

\noindent Moreover, if we consider a reversion $(X, D) \dasharrow (X', D')$ of $1$-standard pairs centered in a point $p \in C_0 \backslash C_1$, we associate to this point its image $(\lambda: 1: 0)$ in $\p^2$ via the map $\tau \circ \eta: X \to \p^2$.

\bigskip
We recall the following lemma, which is an important tool to compute the graph $\F_V$ explicitly:

\bigskip
\begin{lemma}\label{JonquieresAutomorphisms2} (\cite{BD}, Lemma 5.2.1) For $i = 1, 2$, let $(X_i, D_i, \overline{\pi}_i)$ be a $1$-standard pair with a minimal resolution of singularities $\mu_i: (Y_i, D_i, \overline{\pi}_i \circ \mu_i) \to (X_i, D_i, \overline{\pi}_i)$ and let $\eta_i: Y_i \to \f_1$ be the (unique) birational morphism. Then, the following statements are equivalent:
\begin{itemize}
\item[(a)] The $\a^1$-fibered surfaces $(X_1 \backslash D_1, \pi_1)$ and $(X_2 \backslash D_2, \pi_2)$ (respectively the pairs $(X_1, D_1, \overline{\pi}_1)$ and $(X_2, D_2, \overline{\pi}_2)$) are isomorphic.
\item[(b)] There exists an element of $\Jon$ (respectively of $\Jon \cap \Aff$) which sends the points blown-up by $\eta_1$ onto those blown-up by $\eta_2$ and sends the curves contracted by $\mu_1$ onto those contracted by $\mu_2$.
\end{itemize}
\end{lemma}

The rest of the article deals with automorphisms of smooth Gizatullin surfaces. Since these automorphisms extend to birational maps between standard completions of these surfaces, it is natural to study birational maps between standard pairs. In a similar way as for $1$-standard pairs, birational maps between standard pairs can be decomposed into fibered modifications and reversions. Preliminary, we show the following special version of Theorem \ref{FactorizationOfBirationalMaps}:

\begin{lemma}\label{FactorizationOfBirationalMaps2} Let $\varphi: (X, D) \dasharrow (X', D')$ be a birational map between $1$-standard pairs. Then there exists a decomposition

$$\varphi = \varphi_m \circ \cdots \circ \varphi_1: (X, D) = (X_0, D_0) \stackrel{\varphi_1}{\to} (X_1, D_1) \stackrel{\varphi_2}{\to} \cdots \stackrel{\varphi_m}{\to} (X_m, D_m) = (X', D')$$

\noindent such that:
\begin{itemize}
\item[(1)] Each map $\varphi_i$ is either a reversion or a fibered map.
\item[(2)] Each reversion $\varphi_i$ is centred in $\lambda = 0$ as well as its inverse ${\varphi'}^{-1}_i$.
\end{itemize}
\end{lemma}

\begin{bew} First, we recall the following fact. Given an arbitrary smooth $1$-standard pair $(X, D)$, a triangular map $h(w_0, z_0) = (aw_0 + P(z_0), bz_0 + c)$ defines in a natural way a fibered map $h: (X, D) \dasharrow (X', D')$. Indeed, $h$ maps the base point set $M_2 \subseteq C_2 \backslash C_1$ of the feathers with mother component $C_2$ on a certain new set $M' \subseteq C_2 \backslash C_1$. Thus, letting $X_1 \to Q$ and $X'_1 \to Q$, respectively, be the blowups of $M_2$ and $M'_2$, respectively, the lift of $h$ gives a fibered map $(X_1, D_1) \dasharrow (X'_1, D'_1)$. Proceeding in the same way for all points blown up by $X \to Q$, we obtain a new $1$-standard pair $(X', D')$, which we denote by $h.(X, D)$.

In the following we abbreviate $h_a(w_0, z_0) := (w_0 + az_0, z_0)$, $a \in \c$. By Theorem \ref{FactorizationOfBirationalMaps} we decompose $\varphi$ into 

$$\varphi = \varphi_m \circ \cdots \circ \varphi_1: (X, D) = (X_0, D_0) \stackrel{\varphi_1}{\to} (X_1, D_1) \stackrel{\varphi_2}{\to} \cdots \stackrel{\varphi_m}{\to} (X_m, D_m) = (X', D'),$$

\noindent where each $\varphi_i$ is a fibered modification or a reversion. Let $i$ be the smallest index such that $\varphi_i$ is a reversion and let $\varphi_i$ be centred in $\lambda = a$ and let its inverse $\psi_i := \varphi^{-1}_i$ be centred in $\lambda = b$. Write $\varphi_i$ as 

$$\varphi_i = h_a \circ (h_{-a} \circ \varphi_i \circ h_b) \circ h_{-b}.$$

\noindent The map $h_{-a} \circ \varphi_i \circ h_b: h_{-b}.(X_{i - 1}, D_{i - 1}) \dasharrow h_a.(X_i, D_i)$ is a reversion which is centred in $\lambda = 0$ as well as its inverse. Let $\psi_3 := h_a$, $\psi_2 := h_{-a} \circ \varphi_i \circ h_b$ and $\psi_1 := h_{-b}$. Consider the new decomposition

$$\varphi = \varphi_m \circ \cdots \circ \varphi_{i - 1} \circ \psi_3 \circ \psi_2 \circ \psi_1 \circ \varphi_{i + 1} \circ \cdots \circ \varphi_1$$

\noindent and use now induction on the number of reversions which do not satisfy property (2). 
\end{bew}

Let us deduce the decomposition theorem for birational maps of standard pairs. Every smooth standard pair $(X, D)$ arises as a blowup of $\p^1 \times \p^1$ as well as every $1$-standard pair $(X', D')$ arises as a blowup of the Hirzebruch surface $\f_1$. On each of these surfaces we introduced affine charts, namely $\f_1 \backslash (C_0 \cup C_1) \cong \a^2$ and $\p^1 \times \p^1 \backslash (C_0 \cup C_1) \cong \a^2$, endowed with affine coordinates $(w_0, z_0)$ and $(x_0, y_0)$. Now, blowing up an arbitrary point $p = (\lambda: 1: 0) \in C_0 \backslash C_1 \subseteq \f_1$ with exceptional curve $E$ and contracting the proper transform of $C_0$ yields $E^2 = 0$ and $C^2_1 = 0$ on the resulting surface. In other words, we get $Q = \p^1 \times \p^1$.

\bigskip
\begin{tikzpicture}

	\draw (0, 0) rectangle (7, 5) node[anchor=north east]{$\f_1$};
	\begin{scope}[thin]
		\draw (1.5, 0.5) -- (1.5, 4.5) node[anchor=south] {$C_0$};					
		\draw (0.5, 1) -- node[anchor=north] {$C_1$} (6.5, 1);					
		\draw (5.5, 0.5) -- (5.5, 4.5) node[anchor=south] {$L$};					
		\draw[dashed] (1, 2.9) -- node[below=2pt] {$1$} (6, 4.1);					
	\end{scope}	

	\draw [->] (5.3, 3.8) -- (5.3, 2.8) node[anchor=north east] {$w_0$};			
	\draw [->] (5.3, 3.8) -- (4.3, 3.8) node[anchor=south east] {$z_0$};  		
	
	\draw [->] (1.7, 2.8) -- (1.7, 1.8) node[anchor=north west] {$s_0$};			
	\draw [->] (1.7, 2.8) -- (2.6, 3) node[anchor=north west] {$t_0$};  		

	\filldraw [black] (5.5, 4) circle (2pt);
	\filldraw [black] (1.5, 3) circle (2pt) node[anchor=south east] {$\lambda$};

	\draw [->, dashed] (7.2, 2.5) -- (7.8, 2.5);

	\draw (8, 0) rectangle (15, 5);

	\begin{scope}[thin]
		\draw (9.5, 4.5) -- node[anchor=east] {$E$} (8.5, 2);					
		\draw (8.5, 3) -- node[anchor=east] {$C_0$} (9.5, 0.5);					
		\draw (8.5, 1) -- node[anchor=north] {$C_1$} (14.5, 1);					
		\draw (13.5, 4.5) -- (13.5, 0.5) node[anchor=north] {$L$};					
		\draw[dashed] (8.5, 4) -- node[below=2pt] {$0$} (14, 4);					
	\end{scope}	

	\draw [->] (13.3, 3.8) -- (13.3, 2.8) node[anchor=north east] {$w_0$};			
	\draw [->] (13.3, 3.8) -- (12.3, 4.2) node[anchor=south east] {$z_0$};  		
	
	\draw [->] (8.9, 2.5) -- (9.2, 3.2) node[anchor=north west] {$t_1$};			
	\draw [->] (8.9, 2.5) -- (9.2, 1.8) node[anchor=north west] {$s_1$};  		

	\filldraw [black] (13.5, 4) circle (2pt);
	\filldraw [black] (8.7, 2.5) circle (2pt);

\end{tikzpicture}

\bigskip
\begin{tikzpicture}

	\draw [->] (0, 2.5) -- (0.6, 2.5);

	\draw (1, 0) rectangle (8, 5) node[anchor=north east]{$Q$};
	\begin{scope}[thin]
		\draw (2.5, 0.5) -- (2.5, 4.5) node[anchor=south] {$E$};					
		\draw (1.5, 1) -- node[anchor=north] {$C_1$} (7.5, 1);					
		\draw (6.5, 0.5) -- (6.5, 4.5) node[anchor=south] {$C_2$};					
		\draw[dashed] (1.5, 4) -- node[below=2pt] {$0$} (7.5, 4);					
	\end{scope}	

	\draw [->] (6.3, 3.8) -- (6.3, 2.8) node[anchor=north east] {$w_0$};			
	\draw [->] (6.3, 3.8) -- (5.5, 4.2) node[anchor=south east] {$z_0$};  		

	\draw [->, dashed, color=red] (6.4, 3.9) -- (6.4, 2.9) node[anchor=west] {$x_0$};			
	\draw [->, dashed, color=red] (6.4, 3.9) -- (5.4, 3.9) node[anchor=north east] {$y_0$};  		
	
	\draw [->] (2.7, 1.2) -- (3.7, 1.2) node[anchor=south west] {$x'$};			
	\draw [->] (2.7, 1.2) -- (2.7, 2.2) node[anchor=south west] {$y'$};  		

	\filldraw [black] (6.5, 4) circle (2pt);
	\filldraw [black] (2.5, 1) circle (2pt);

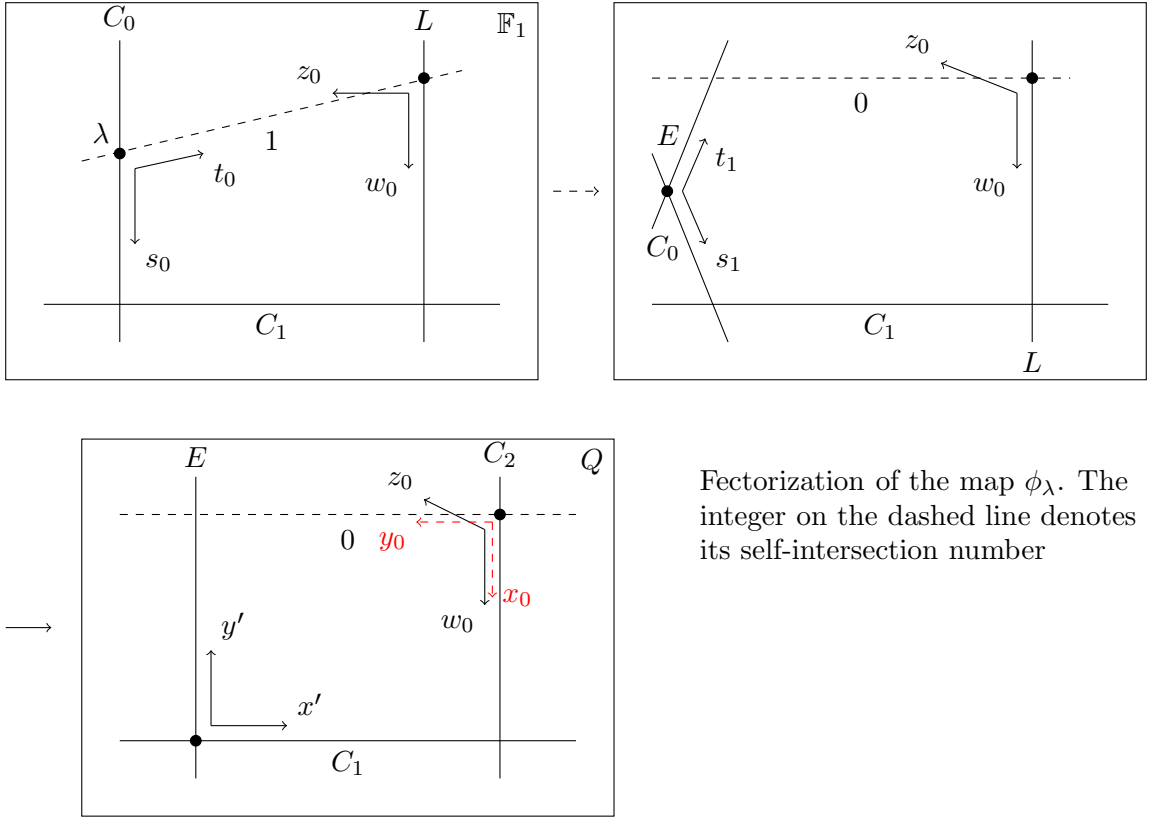
\node[draw=white, align=left] at (12,4) {Fectorization of the map $\phi_\lambda$. The \\ integer on the dashed line denotes \\ its self-intersection number};

\end{tikzpicture}

\bigskip
Let us denote this map by $\phi_\lambda: (\f_1, C_0 \cup C_1) \dasharrow (\p^1 \times \p^1, C_0 \cup C_1)$. However, performing these transformations we obtain an affine coordinate system $(x', y')$, centred in $C_0 \cap C_1$ (cf. the figure above). Introducing the coordinates $(s_0, t_0) = (w_0/z_0 - \lambda, 1/z_0)$, $(s_1, t_1) = (s_0, t_0/s_0)$ and $(x', y') = (s_1t_1, t_1)$, we obtain

$$(x', y') = \left( \frac{1}{z_0}, \frac{1}{w_0 - \lambda z_0} \right).$$

\noindent In other words, we get $(x_0, y_0) = \left( 1/y', 1/x' \right) = (w_0 - \lambda z_0, z_0)$, up to an isomorphism $(x_0, y_0) \mapsto (\lambda x_0, \mu y_0)$ of $\p^1 \times \p^1$. In particular, this implies the following. Given finite subsets $M_3, \dots, M_{n - 1} \subset \c^*$, we can construct a standard pair $(X, D)$ as in \ref{CoordinatesGeneralCase} as well as a $1$-standard pair $(X', D')$ applying precisely the same algorithm, but using the coordinates $(w_0, z_0)$ on $\f_1$ instead of $(x_0, y_0)$ on $\p^1 \times \p^1$ (this yields the same Gizatullin surface $V = X \backslash D = X' \backslash D'$). Then $(X^\vee, D^\vee) = (\phi^{-1}_0 \circ \varphi \circ \phi_0)(X', D')$, where $\varphi$ denotes the reversion of a standard pair. 

Now, similar to the case of $1$-standard pairs we say that a birational map $\varphi: (X, D) \dasharrow (X', D')$ between standard pairs is a \emph{fibered modification}, if it is a lift of a triangular map $\p^1 \times \p^1 \backslash (C_0 \cup C_1) \to \p^1 \times \p^1 \backslash (C_0 \cup C_1)$, $(x_0, y_0) \mapsto (ax_0 + P(y_0), by_0 + c)$, which is not an isomorphism (note, that in the case of standard pairs such maps give isomorphisms between the cooresponding completions if and only if $\deg(P) = 0$). Hence, Lemma \ref{FactorizationOfBirationalMaps2} and the considerations above immediately imply the following corollary on the decomposition of birational maps between standard pairs:

\begin{cor}\label{FactorizationOfBirationalMapsStandard} Let $\varphi: (X, D) \dasharrow (X', D')$ be a birational map between standard pairs. Then there exists a decomposition

$$\varphi = \varphi_m \circ \cdots \circ \varphi_1: (X, D) = (X_0, D_0) \stackrel{\varphi_1}{\to} (X_1, D_1) \stackrel{\varphi_2}{\to} \cdots \stackrel{\varphi_m}{\to} (X_m, D_m) = (X', D')$$

\noindent such that each $\varphi_i$ is either a reversion or a fibered modification.
\end{cor}

\subsection{Smooth Gizatullin surfaces and automorphisms of $\a^1$-fibrations}

\bigskip
There is a good chance to obtain $\a^1$-fibered affine surfaces which behaves well under applying birational maps, if the extended divisor admits as many as possible inner components. For the rest of the article we will assume that $V$ is a smooth Gizatullin surface satisfying the following condition:

\begin{eqnarray*}
(*) && V \ \text{admits a} \ \text{standard completion} \ (X, D) \ \text{with} \ n \geq 4 \ \text{such that} \ C_3, \dots, C_{n - 1} \ \text{are} \\
&& *\text{-components and there is no feather attached to} \ C_2 \ \text{and to} \ C_n.
\end{eqnarray*}

Using the Matching Principle it is not hard to see that if $V$ satisfies $(*)$, then \emph{every} standard completion of $V$ also satisfies this condition. The same holds for any $m$-standard completion of $V$, since $m$-standard completions of $V$ are obtained from standard completions by performing elementary transformations centred in $C_0$.\\

\begin{thm}\label{SurfacesRigidExtendedDivisor} Let $V$ be as in $(*)$. Then the following hold:
\begin{itemize}
\item[(1)] For any two $1$-standard completions $(X', D')$ and $(X'', D'')$ of $V$ with $(X' \backslash D', \pi') \cong (X'' \backslash D'', \pi'')$ we have $(X', D', \bar{\pi}') \cong (X'', D'', \bar{\pi}'')$. 
\item[(2)] The graph $\F_V$ has one of the following two forms:

$$\F_V: \begin{xy}
  \xymatrix{
  [(X, D)] \ \bullet \ar@{<->}[r] & \bullet \ [(X^\vee, D^\vee)] \\
  }
\end{xy} \quad \text{or} \quad \F_V: [(X, D)] \ \bullet \rcirclearrowleft .$$

\noindent If $\F_V$ is of the form $\bullet \rcirclearrowleft$, then $D^{\geq 2}$ is a palindrome and $r_i = r_{i^\vee}$ and $Q_i = Q_{i^\vee}$ holds for all $i = 3, \dots, n - 1$.
\item[(3)] If $r_i > 0$ holds for at most two indices $i \in \{ 3, \dots, n - 1 \}$, then $\F_V$ has the form $\bullet \rcirclearrowleft$ if and only if $D^{\geq 2}$ is a palindrome and $r_i = r_{i^\vee}$ and $Q_i = Q_{i^\vee}$ holds for all $i = 3, \dots, n - 1$.
\item[(4)] $\Aut(V)$ is generated by automorphisms of $\a^1$-fibrations if and only if $\F_V$ has no loops except for the case $\Gamma_D = [[0,-1, -2, -2, -2]]$.
\end{itemize}
\end{thm}

\begin{bew} We consider two $1$-standard completions $(X', D')$ and $(X'', D'')$ of $V$ such that $(X' \backslash D', \pi') \cong (X'' \backslash D'', \pi'')$. By Lemma \ref{JonquieresAutomorphisms2} an isomorphism $\Jon \ni \psi: (X' \backslash D', \pi') \stackrel{\sim}{\to} (X'' \backslash D'', \pi'')$ sends the centers of the blowup $\eta': X' \to \f_1$ onto those of $\eta'': X'' \to \f_1$. According to Lemma \ref{JonquieresAutomorphisms} we present $\psi$ in the coordinates $(w_0, z_0)$ of $\a^2 = \f_1 \backslash (C_0 \cup C_1)$ as

$$\psi(w_0, z_0) = (aw_0 + P(z_0), bz_0 + c), \quad a, b \in \c^*, c \in \c, P(z_0) \in \c[z_0].$$

Assuming that the blowup process for both surfaces starts at $(0, 0) \in L_0$, we obtain $c = P(0) = 0$ and $\psi$ has the form

$$\psi(w_0, z_0) = (aw_0 + z_0Q(z_0), bz_0), \quad a, b \in \c^*, \ Q(y_0) \in \c[z_0].$$

We decompose $\eta'$ into single blowups $X' = X'_{n - 2 + m} \to \cdots \to X'_{n - 2} \to \cdots \to X'_0 = \f_1$, such that we create the boundary components by inner blowups in $X'_{n - 2} \to \cdots \to X'_0 = \f_1$ and the $m = r_3 + \cdots + r_{n - 1}$ feathers in $X'_{n - 2 + m} \to \cdots \to X'_{n - 2}$ (similarly for $X''$). We fix an index $j \in \{ 3, \dots, n - 1 \}$ and consider $\psi$ on the component $C_j$ on the surface $X_{n - 2}$. In each step we introduce inductively affine coordinates $(w_i, z_i)$ on $X_i$, such that:

\begin{itemize}
\item[(1)] $C_j \backslash C_{j - 1} = \{ z_{n - 2} = 0 \}$ and $C_{j + 1} \backslash C_{j + 2} = \{ w_{n - 2} = 0 \}$ on $X_{n - 2}$ (here we define $C_{j + 2}$ to be the infinity point if $j = n - 1$),
\item[(2)] after each blowup we have either $(w_i, z_i) = (w_{i + 1}, z_{i + 1})$ or $(w_i, z_i) = (w_{i + 1}, w_{i + 1}z_{i + 1})$ or $(w_i, z_i) = (w_{i + 1}z_{i + 1}, z_{i + 1})$.
\end{itemize}

We refer to the last two coordinate transformations as to transformations of \emph{type 1} and \emph{type 2} respectively. Denoting the lift of $\psi$ on the surface $X_i$ by $\Psi_i$, we show by induction that $\Psi_{n - 2}$ has the form

$$\Psi_{n - 2}(w_{n - 2}, z_{n - 2}) = \left( \alpha w_{n - 2}(1 + w_{n - 2}^kz_{n - 2}^lR(w_{n - 2}^sz_{n - 2}^t))^p, \frac{\beta z_{n - 2}}{(1 + w_{n - 2}^kz_{n - 2}^lR(w_{n - 2}^sz_{n - 2}^t))^q} \right),$$

\noindent such that $R = \gamma \cdot Q$ with an appropriate $\gamma \in \c^*$, $k \geq 0$, $l, s, t, p, q \geq 1$ and $\alpha$ and $\beta$ are monomials in $a$ and $b$. In the first step we blow up $\f_1$ in the point $(0, 0) \in L_0$ and introduce the coordinates $(w_1, z_1)$ via $(w_0, z_0) = (w_1, w_1z_1)$. This leads to

$$\Psi_1(w_1, z_1) = \left( w_0(a + z_0Q(w_0z_0)), \frac{bz_0}{a + z_0Q(w_0z_0)} \right) = \left( aw_0(1 + z_0R(w_0z_0)), \frac{(b/a)z_0}{1 + z_0R(w_0z_0)} \right)$$

\noindent with $R = \frac{1}{a}Q$. In the induction step we consider the case $(w_i, z_i) = (w_{i + 1}, w_{i + 1}z_{i + 1})$. A short computation yields

$$\Psi_{i + 1}(w_{i + 1}, z_{i + 1}) = \left( \alpha w_{i + 1}(1 + w_{i + 1}^{k + l}z_{i + 1}^lR(w_{i + 1}^{s + t}z_{i + 1}^t))^p, \frac{(\beta/\alpha) z_{i + 1}}{(1 + w_{i + 1}^{k + l}z_{i + 1}^lR(w_{i + 1}^{s + t}z_{i + 1}^t))^{p + q}} \right).$$

\noindent Similarly we obtain in the case $(w_i, z_i) = (w_{i + 1}z_{i + 1}, z_{i + 1})$

$$\Psi_{i + 1}(w_{i + 1}, z_{i + 1}) = \left( (\alpha/\beta) w_{i + 1}(1 + w_{i + 1}^kz_{i + 1}^{k + l}R(w_{i + 1}^sz_{i + 1}^{s + t}))^{p + q}, \frac{\beta z_{i + 1}}{(1 + w_{i + 1}^kz_{i + 1}^{k + l}R(w_{i + 1}^sz_{i + 1}^{s + t}))^q} \right).$$

The case $(w_i, z_i) = (w_{i + 1}, z_{i + 1})$ is obvious (we need such transformations only when we perform inner blowups centred outside the affine piece with coordinates $(w_i, z_i)$). Moreover, we see from the induction step that $k = 0$ holds if and only if we perform no blowups of type 1 except for the first one in $(0, 0) \in L_0$.

Hence, $\Psi_{n - 2}$ induces on $C_j \backslash (C_{j - 1} \cup C_{j + 1}) \cong \c^*$ the map $\Psi_{n - 2}(w, 0) = (\alpha w, 0)$. But the affine map 

$$\tilde{\psi}(w_0, z_0) := (aw_0 + cz_0, bz_0) \quad \text{with} \quad c = Q(0),$$

\noindent defines the same map on $C_j \backslash (C_{j - 1} \cup C_{j + 1})$. In particular, $\tilde{\psi}$ sends the points blown-up by $\eta': X' \to \f_1$ onto those blown-up by $\eta'': X'' \to \f_1$. By Lemma \ref{JonquieresAutomorphisms2} we obtain that $(X', D', \bar{\pi}') \cong (X'', D'', \bar{\pi}'')$. This shows (1).

We show that the graph $\F_V$ admits only one arrow (the first part of assertion (2) follows immediately since $\F_V$ is connected). Due to (1) we can choose $(X, D, \bar{\pi})$ itself as a representative of the conjugacy class $[(X, D)] \in \F_V$. The automorphism

$$\psi_a(w_0, z_0) = (w_0 + az_0, z_0), \quad a \in \c,$$

\noindent of $\a^2 = \f_1 \backslash (C_0 \cup C_1)$ can be lifted to an automorphism of $(X, D)$. Moreover, $\psi_a$ induces on $C_0 \backslash C_1 \cong \a^1$ the translation $\lambda \mapsto \lambda + a$. Thus $\Aut(X, D)$ acts transitively on $C_0 \backslash C_1$ and there is only one arrow starting from $[(X, D)]$.

The condition on $\F_V$ to have a loop is equivalent to $(X, D, \bar{\pi}) \cong (X^\vee, D^\vee, \bar{\pi}^\vee)$. This implies the second part of (2).

To show assertion (3), we consider two $1$-standard completions $(X', D')$ and $(X'', D'')$ of $V$ such that there exist two indices $s, t \in \{ 3, \dots, n - 1 \}$, $s < t$, with  $r'_i, r''_i > 0$ only for $i \in \{ s, t \}$ and $(Q'_s, Q'_t) = (Q''_s, Q''_t)$. It is sufficient to show that this implies $(X', D', \bar{\pi}') \cong (X'', D'', \bar{\pi}'')$. Then we get the non-trivial direction of (3) as follows. If $(X, D, \bar{\pi})$ is a $1$-standard completion of $V$ with $r_i = r_{i^\vee} > 0$, $r_j = 0$ for $j \neq i, i^\vee$ and $Q_i = Q_{i^\vee}$, the Matching Principle yields $r^\vee_{i^\vee} = r_i$, $r^\vee_i = r_{i^\vee}$, $Q^\vee_{i^\vee} = Q_i$ and $Q^\vee_i = Q_{i^\vee}$ for the corresponding invariants. Thus we obtain $r^\vee_k = r_k$ and $Q^\vee_k = Q_k$, $k \in \{ i, i^\vee \}$ and $(X, D, \bar{\pi}) \cong (X^\vee, D^\vee, \bar{\pi}^\vee)$ follows. Hence the graph $\F_V$ admits only one vertex.

Again, we consider the blowup process $(X'_{n - 2 + m}, D'_{n - 2 + m}) \to \cdots \to (X'_{n - 2}, D'_{n - 2}) \to \cdots \to (X'_0, D'_0) = (\f_1, C_0 \cup C_1)$, where we create the boundary components in $(X'_{n - 2}, D'_{n - 2}) \to \cdots \to (X'_0, D'_0) = (\f_1, C_0 \cup C_1)$ and the $m = r'_s + r'_t$ feathers in $(X'_{n - 2 + m}, D'_{n - 2 + m}) \to \cdots \to (X'_{n - 2}, D'_{n - 2})$. The $2$-torus action on $\f_1 \backslash (C_0 \cup C_1)$ lifts to $(X'_{n - 2}, D'_{n - 2})$, inducing maps $x \mapsto \alpha_i x$ on $C_i \backslash (C_{i - 1} \cup C_{i + 1}) \cong \c^*$ (we identify $C_i \backslash (C_{i - 1} \cup C_{i + 1})$ with $\c^*$ in such a way that $C_i \cap C_{i - 1} = \{ \infty \}$ and $C_i \cap C_{i + 1} = \{ 0 \}$). Since no feathers are attached to $C_i$ for $i \neq s, t$, these maps do not play any role on these components. Further, the maps $x \mapsto \alpha_s x$ and $x \mapsto \alpha_t x$ can be chosen arbitrarily. Indeed, it can be checked by induction on $n$ that an element $(a, b)$ of the $2$-torus $\t := \{ (w_0, z_0) \mapsto (aw_0, bz_0) \mid a, b \in \c^* \} \subseteq \Aff \cap \Jon$, which lifts to an isomorphism of pairs $(X', D') \stackrel{\sim}{\to} (X'', D'')$, induces on $X'_{n - 2}$ maps $C'_j \backslash (C'_{j - 1} \cup C'_{j + 1}) \to C''_j \backslash (C''_{j - 1} \cup C''_{j + 1})$ of the form 

\begin{eqnarray}\label{ArbitraryMultiples}
w'_{n - 2} \mapsto w''_{n - 2} = a^{p_i}b^{q_i}w'_{n - 2} \quad \text{such that} \quad p_iq_j - p_jq_i \neq 0
\end{eqnarray}

\noindent for $i \neq j$. Hence there exist suitable $a, b \in \c^*$ such that $a^{p_s}b^{q_s} = \alpha_s$ and $a^{p_t}b^{q_t} = \alpha_t$. Thus letting $A'_s, A'_t, A''_s$ and $A''_t$ be the corresponding base point sets of the feathers, there exists a $2$-torus action $\mu$ such that $\mu.A'_s = A''_s$ and $\mu.A'_t = A''_t$. Therefore we get $(X', D') \cong (X'', D'')$.

To show (4), we let $\Gamma_D \neq [[0, -1, -2, -2, -2]]$. It is easy to see that this condition is equivalent to the existence of a boundary component of self-intersection $\leq -3$. In this case assertion (4) is a direct consequence of Proposition \ref{SurfaceTree} and assertion (2).

It remains the case where $\Gamma_D = [[0, -1, -2, -2, -2]]$. We choose a fixed reversion $\Psi_0: (X, D) \dasharrow (X, D)$ with center $p \in C_0 \backslash C_1$ and such that $p$ is also the center of $\Psi_0^{-1}$ (we can choose such a reversion since $\Aut(X, D)$ acts transitively on $C_0 \backslash C_1$). Let $\alpha \in \Aut(X, D)$ be an element, which does not fix the center of $\Psi_0$. Then the reversions $\Psi^{-1}_0 = \Psi_0$ and $\Psi_0 \alpha$ have distinct base points, \ie $\Psi_0 \alpha \Psi_0$ is a reversion by Lemma \ref{CompositionOfReversions}, equal to $\beta \Psi_0 \gamma$ for some $\beta, \gamma \in \Aut(X, D)$. Therefore we have $\Psi_0 = (\alpha^{-1}\Psi_0)(\beta \Psi_0 \gamma) = \alpha^{-1}(\Psi_0 \beta \Psi^{-1}_0) \gamma$. Since $\beta$ preserves the fibration $\pi$, there is a $\varphi \in \Aut(\a^1)$ such that $\pi \beta = \varphi \pi$. This yields $(\pi \Psi_0)(\Psi_0 \beta \Psi^{-1}_0) = \varphi(\pi \Psi_0)$, \ie the map $\Psi_0 \beta \Psi^{-1}_0$ is compatible with the fibration $\pi \Psi_0$. Thus the reversion $\Psi_0$ is generated by automorphisms of $\a^1$-fibrations. In Corollary \ref{StructureOfAutomorphismGroup} below we show that $\Aut(V)$ is generated by $\langle \Aut(X, D), \Psi_0 \rangle$ and $\Aut(V, \pi)$ (the elements of $\Aut(X, D)$ as well as $\Psi_0$ give automorphisms of $V$ by restriction). It follows that $\Aut(V)$ is generated by automorphisms of $\a^1$-fibrations.
\end{bew}

\begin{rem} (1) Theorem \ref{SurfacesRigidExtendedDivisor} yields in particular that every $\a^1$-fibration $\varphi: V \to \a^1$ on such a surface $V$ is conjugated either to $\Phi_0 := \Phi_{|C_0|}: V \to \a^1$ or to $\Phi^\vee_0 := \Phi_{|C^\vee_0|}: V \to \a^1$. This is a special case of Theorem 5.10 in \cite{FKZ3}, where the same assertion is shown for normal Gizatullin surfaces with a so called \emph{distinguished} and \emph{rigid} extended divisor. In the terminology of \cite{FKZ3} these are, considering the smooth case, precisely those Gizatullin surfaces, which satisfy $(*)$.

(2) Conversely, the condition that $V$ admits at most two conjugacy classes of $\a^1$-fibrations  does \textit{not} imply that $V$ admits an extended divisor which satisfies $(*)$. We consider the following example. Let $V$ be a smooth Gizatullin surface admitting a $1$-standard completion with extended divisor

\bigskip
$$D_{\text{ext}}: \quad \cu{0} \lin \cu{-1} \lin \cu{-2} \lin \cu{-2} \nlin \cshiftup{-1}{} \lin \cu{-4} \lin \cu{-2} \nlin \cshiftup{-1}{} \lin \cu{-2} \quad .$$

Note, that the components $C_2, C_4$ and $C_6$ are $+$-components, while $C_3$ and $C_5$ are of type $*$. We can construct $(X, D)$ from $(\f_1, C_0 \cup C_1)$ in the following way: we blow up two times in $(0, 0) \in L_0$ and introduce coordinates $(w_1, z_1)$ via $(w_0, z_0) = (w_1, w_1^2z_1)$ as well as coordinates $(u_1, v_1)$ via $(w_0, z_0) = (u_1v_1, v_1)$. This leads to the dual graph

$$\cu{0} \lin \cu{-1} \lin \cu{-2} \lin \cu{-1} \lin \cu{-2} \quad .$$

Now we perform an outer blowup in a point $P = (\beta, 0) \in C_4 \backslash C_3$ (in coordinates $(u_1, v_1)$), obtaining the component $C_6$. Finally we perform an inner blowup in $C_4 \cap C_6$, which results in the component $C_5$. After the blowup in $P$ we introduce the coordinates $(u_2, v_2)$ via $(u_1 - \beta, v_1) = (u_2, u_2v_2)$ and after the last blowup the coordinates $(u_3, v_3)$ via $(u_2, v_2) = (u_3, u_3v_3)$. We denote the resulting surface by $(\tilde{X}, \tilde{D})$. In the last step we create both $(-1)$-feathers by blowing up in some points $Q_1 = (0, \alpha) \in C_3 \backslash (C_2 \cup C_4)$ (in coordinates $(w_1, z_1)$) and $Q_2 = (0, \gamma) \in C_5 \backslash (C_4 \cup C_6)$ (in coordinates $(u_3, v_3)$). An appropriate automorphism $\psi(w_0, z_0) = (aw_0 + bz_0, cz_0)$ of $(\f_1, C_0 \cup C_1)$ can bring the surface $(X, D)$ in a "standard form": The condition $b = -a\beta$ moves the point $(\beta, 0)$ on $(0, 0)$ (in coordinates $(u_1, v_1)$). Thus we may assume that $\beta = 0$. Then $\psi$ can be lifted to $(\tilde{X}, \tilde{D})$ if and only if $b = 0$ holds and the lifting $\Psi$ of $\psi$ has the following forms in the coordinates introduced above:

$$\Psi(w_1, z_1) = (aw_1, a^{-2}cz_1) \quad \text{and} \quad \Psi(u_3, v_3) = (ac^{-1}u_3, a^{-2}c^3)$$

The conditions $\Psi(Q_1) = (0, 1)$ and $\Psi(Q_2) = (0, 1)$ (in the corresponding affine coordinates) leads to $c = a^2$ and $c^2 = 1$, or equivalently $(a, c) \in \{ (1, 1), (-1, 1), (i, -1), (-i, -1) \}$. In particular, $ac^{-1}$ can take every value in $W_4 = \{ z \in \c^* \mid z^4 = 1 \} = \{ \pm 1, \pm \text{i} \}$.
Therefore we obtain: 

\begin{itemize}
\item[(a)] $V$ does not depend on any parameter.
\item[(b)] $(X, D) \cong (X^\vee, D^\vee)$.
\item[(c)] For any two $1$-standard completions $(X, D)$ and $(X', D')$ of $V$ it follows $(X, D) \cong (X', D')$ (see Lemma \ref{JonquieresAutomorphisms2}). In particular, $V$ admits only one conjugacy class of $\a^1$-fibrations.
\end{itemize}

The subgroup $\Aut(X, D)$ does \emph{not} act transitively on $C_0 \backslash C_1$. If we identify $C_0 \backslash C_1$ with $\a^1$ via $(\lambda: 1: 0) \mapsto \lambda$, then the orbit of a point $\lambda \in C_0 \backslash C_1$ is given by $\Aut(X, D).\lambda = \lambda \cdot W_4$ ($\psi(w_0, z_0) = (aw_0, cz_0)$ induces on $C_0 \backslash C_1$ the map $\lambda \mapsto ac^{-1} \cdot \lambda$). Since there is a $1:1$-correspondence betweens arrows and $\Aut(X, D)$-orbits of $C_0 \backslash C_1$ (see Remark \ref{EquivalenceOfArrows}), the graph $\F_V$ has the form

$$\F_V: \quad \lambda_1 \rcirclearrowright \stackrel{\cdots}{\bullet} \lcirclearrowleft \lambda_2 \quad ,$$

\noindent where the arrows are in a $1:1$-correspondence to elements of $\c/W_4$.
\end{rem}

\begin{cor} Let $V$ be as in Theorem \ref{SurfacesRigidExtendedDivisor} and let $n$ be odd. Then $\Aut(V)$ is generated by automorphisms of $\a^1$-fibrations.
\end{cor}

\begin{bew} We claim that if the subgraph $D_{\text{ext}}^{\geq 2}$ of an extended divisor $D_{\text{ext}}$ is symmetric, then $n$ is even. 

We show by induction on $n$ that this cannot occur for odd $n$. In the case $n = 5$ we can blow up the boundary divisor of $\f_1$ (which is of the type $[[0, -1, 0]]$) either to a zigzag of the type $[[0, -1, -3, -1, -2, -2]]$ or of the type $[[0, -1, -3, -1, -2, -2]]$. However, further blowups cannot produce a symmetric extended divisor.\\
Assume, that there is an odd number $n > 5$ such that $D_{\text{ext}}^{\geq 2}$ is symmetric. Consider the zigzag $D'$, which arises after blowing down all feathers of $D_{\text{ext}}$. The dual graph $\Gamma_{D'}$ of $D'^{\geq 2}$ is also symmetric, say of the form

$$\Gamma_{D'} = [[0, -1, w_2, \dots, w_k, w_k, \dots, w_2]], \quad n = 2k - 1.$$

Since $D_{\text{ext}}^{> 2}$ is contractible, one of the $w_i$ is equal to $-1$. If $w_i = -1$ for some $i \in \{ 2, \dots, k - 1 \}$, then we can blow down the corresponding $(-1)$-curves $C_i$ and $C_{i^\vee}$ and obtain a symmetric boundary divisor $\tilde{D} = \tilde{C}_0 \cup \cdots \cup \tilde{C}_{n - 2}$ of lenght $n - 2$. By generating symmetrically feathers on $\tilde{D}$ (until $\tilde{C}^2_i \leq -2$ for all $2 \leq i \leq n - 2$) we obtain an extended divisor $\tilde{D}_{\text{ext}}$, such that $\tilde{D}_{\text{ext}}^{\geq 2}$ is symmetric. This is impossible by induction hypothesis. The case $w_k = -1$ cannot occur, since otherwise $\Gamma_{D'}$ cannot be blown down to $[[0, -1, 0]]$.
\end{bew}

\subsection{Invariant subsets under the action of $\Aut(V)$}

Our next goal is to describe the action of the automorphism group of smooth Gizatullin surfaces as in Theorem \ref{SurfacesRigidExtendedDivisor}. We will show that these surfaces in general admit points, which do not belong to the big orbit $O$ of $\Aut(V)$. These points are even fixed points of $\Aut(V)$, if the base points of the feathers are "in general position". We start with the following simple observation:

\begin{prop}\label{AutFiniteComplement} Let $V$ be a smooth Gizatullin surface with a standard completion $(X, D)$ and the associated $\p^1$-fibration $\overline{\pi} := \Phi_0: X \to \p^1$. Furthermore, let $F_i$ be the feathers of the extended divisor $D_{\text{ext}}$ and let $F^\vee_i$ be the matching feathers in $D^\vee_{\text{ext}}$, where $(X, D) \dasharrow (X^\vee, D^\vee)$ is the reversion. Letting $O$ be the big orbit of the natural action of $\Aut(V)$, we have

$$V \backslash O \subseteq \bigcup_{i, j} (F_i \cap F^\vee_j).$$

\end{prop}

\begin{bew} Since $\Aut(V)$ acts on $V$ with a big orbit $O$, each fiber $F'$ of $\pi = \overline{\pi}\vert_V$ intersects $O$. Since $F_0 := \overline{\pi}^{-1}(\overline{\pi}(C_2 \cup \cdots \cup C_n))$ is the only degenerated fiber of $\overline{\pi}$, every point of $V \backslash O$ cannot belong to any other fiber $F$. Indeed, if $D = C_0 \cup \cdots \cup C_n$, then the automorphism $\psi(x_0, y_0) = (x_0 + y_0^{n - 1}, y_0)$ of $\a^2 = \p^1 \times \p^1 \backslash (C_0 \cup C_1)$ can be lifted to an automorphism of $V$, which induces a non-trivial translation $x \mapsto x + y^{n - 1}$ on any other fiber $F = \pi^{-1}(y) \cong \a^1$, $y \neq 0$. Thus $V \backslash O \subseteq F_0 \cap V$.\\
Let now $x \in V$ be a point contained in

$$\left( \bigcup_i F_i \right) \backslash \left( \bigcup_{i, j} (F_i \cap F^\vee_j) \right).$$

\noindent Then $x$ is contained in some regular fiber of the $\a^1$-fibration $\pi^\vee$, thus it can be moved continuously by an appropriate $\c_+$-action. Therefore $x \in O$. In a similar way we obtain that any point in

$$\left( \bigcup_i F^\vee_i \right) \backslash \left( \bigcup_{i, j} (F_i \cap F^\vee_j) \right)$$

\noindent can be moved to $O$ by automophisms of $V$. It follows that $V \backslash O \subseteq \bigcup_{i, j} (F_i \cap F^\vee_j)$.
\end{bew}

The main result of this paper is a class of smooth Gizatullin surfaces which yield counter-examples to Gizatullin's conjecture. Indeed, we show that the automorphism group of surfaces considered in Theorem \ref{SurfacesRigidExtendedDivisor} does not act transitively in general. First of all, we need the following simple lemma:

\begin{lemma}\label{FiniteSet} Let $A \subseteq \c^*$ be a non-empty finite subset and let $W_m := \{ z \in \c^* \mid z^m = 1 \} = \{ e^{\frac{2k\pi i}{m}} \mid 0 \leq k \leq m - 1 \}$, $m \geq 1$, be the group of $m$-th roots of unity. We represent $A$ as

\begin{equation}\label{ConditionForInterchange}
A = \bigcup_{i = 1}^s c_iW_{m_i}
\end{equation}

with $c_i \in \c^*$, $m_i \geq 1$ and such that $s$ is \emph{minimal}. 
\begin{itemize}
\item[(1)] For the group $G := \{ \alpha \in \c^* \mid \alpha \cdot A = A \}$ it holds 

$$G = W_d \quad \text{with} \quad d = \gcd(m_1, \dots, m_s).$$

\noindent In particular, $G \cong \z_d$ and there exists an $\alpha \in \c^* \backslash \{ 1 \}$ with $\alpha \cdot A = A$ if and only if $d \geq 2$.
\item[(2)] The $G$-action on $A$ yields a decomposition of $A$ in exactly $m(A) := \frac{m_1 + \cdots + m_s}{d}$ orbits $B_1, \dots, B_{m(A)}$.
\end{itemize}
\end{lemma}

We will use the following notation: if we represent a finite non-empty subset $A \subseteq \c^*$ in the form $A = \bigcup_{i = 1}^s c_iW_{m_i}$ such that $s$ is minimal, as in (\ref{ConditionForInterchange}), we let

$$d(A) := \gcd(m_1, \dots, m_s), \quad G(A) := \{ \alpha \in \c^* \mid \alpha \cdot A = A \} \cong \z_{d(A)} \quad \text{and} \quad m(A) := \frac{m_1 + \cdots + m_s}{d(A)}.$$

\noindent In addition, if $A = \emptyset$, we let

$$d(A) = d(\emptyset) := 0, \quad G(A) = G(\emptyset) := \c^* \quad \text{and} \quad m(A) = m(\emptyset) := 0.$$

In order to formulate the main result, we need to determine the boundary components with the property that for every feather $F$ attached to this component we have $F \backslash D \subseteq O$.

\begin{df}\label{DfExceptionalComponent} Let $V$ be as in $(*)$ and let $(X, D)$ be a standard completion of $V$. Then $(X, D)$ arises as a blowup of the quadric $Q = \p^1 \times \p^1$ and a suitable ordering of the blowups yields an intermediate surface $(X', D')$, $D' = C'_0 \cup \cdots \cup C'_{m + 3}$, with boundary divisor of the type $[[0, 0, -2, -2, \dots, -2, -1, -(m + 1)]]$, such that $m$ is maximal. The proper transform of $C'_3, \dots, C'_{m + 2}$ under the blowup $X \to X'$ are boundary components $C_{\tau_1}, \dots, C_{\tau_m}$ of $X$ with $3 \leq \tau_1 < \cdots < \tau_m \leq n - 1$. We call these $C_{\tau_i}$ the \emph{exceptional components of} $D$ and we denote $\mathfrak{E}_D := \{ \tau_1, \dots, \tau_m \}$. Further, we denote the set $\{ \tau_1^\vee, \dots, \tau_m^\vee \}$ by $\mathfrak{E}^\vee_D$. Note, that these $C_{\tau_i}$ are uniquely determined by this condition.
\end{df}

For example, if $D_{\text{ext}}$ has the dual graph

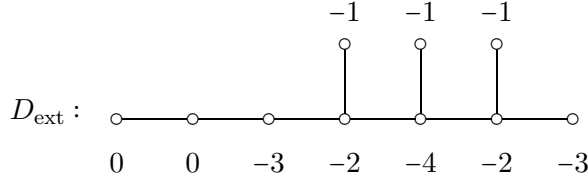
\vspace{15pt}
$$D_{\text{ext}}: \quad \cu{0} \lin \cu{0} \lin \cu{-3} \lin \cu{-2} \nlin \cshiftup{-1}{} \lin \cu{-4} \nlin \cshiftup{-1}{} \lin \cu{-2} \nlin \cshiftup{-1}{} \lin \cu{-3}$$

\noindent then $C_4$ and $C_5$ are the exceptional components of $D$. Indeed, the blowup $X \to Q$ factorizes as follows:

\vspace{15pt}
$$\cou{C_0}{0} \lin \cou{C_1}{0} \lin \cou{C_2}{0} \quad \leftarrow \quad \cou{C_0}{0} \lin \cou{C_1}{0} \lin \cou{C_2}{-2} \lin \cou{C_4}{-2} \lin \cou{C_5}{-1} \lin \cou{C_6}{-3} \quad \leftarrow \quad \cu{0} \lin \cu{0} \lin \cu{-3} \lin \cu{-2} \nlin \cshiftup{-1}{} \lin \cu{-4} \nlin \cshiftup{-1}{} \lin \cu{-2} \nlin \cshiftup{-1}{} \lin \cu{-3}$$

\begin{rem} Note, that if $D^{\geq 2}_{\text{ext}}$ is symmetric, then $C_{\lfloor \frac{n}{2} \rfloor + 1}$ is always exceptional. Moreover, due to the maximality of $m$ it follows that $C_{n - 1}$ is an exceptional component for an arbitrary extended divisor $D_{\text{ext}}$ as well. Note also, that in the case where $D_{\text{ext}}^{\geq 2}$ is symmetric, we have $\mathfrak{E}_D = \mathfrak{E}_{D^\vee}$. 
\end{rem}

\bigskip
The following Theorem is the main result  of this paper:

\bigskip
\begin{thm}\label{FixedPointOfSurface} Let $V$ be as in $(*)$, $(X, D)$ a standard completion of $V$, and let\\
$\bullet$ $A_i = \{ P_{i, 1}, \dots, P_{i, r_i} \} \subseteq C_i \backslash (C_{i - 1} \cup C_{i + 1}) \cong \c^*$, $3 \leq i \leq n - 1$, be the base point set of the feathers $F_{i, j}$,\\ 
$\bullet$ $B_{i, 1}, \dots, B_{i, m(A_i)}$ be the orbits of the $G(A_i)$-action on $A_i$ for $i \not\in \mathfrak{E}_D \cup \mathfrak{E}^\vee_{D^\vee}$.\\
\noindent If $Q(X, D) \neq Q(X^\vee, D^\vee)$, then for $i \not\in \mathfrak{E}_D \cup \mathfrak{E}^\vee_{D^\vee}$, $1 \leq j \leq m(A_i)$, we let 

$$O_{i, j} := \bigcup_{l \in \{ 1, \dots, r_i \}; P_{i, l} \in B_{i, j}} (F_{i, l} \cap F^\vee_{i, l}) \subseteq V,$$

\noindent Otherwise, if $Q(X, D) = Q(X^\vee, D^\vee)$, then for $i \not\in \mathfrak{E}_D \cup \mathfrak{E}^\vee_{D^\vee}$, $i < \lfloor \frac{n}{2} \rfloor + 1, 1 \leq j \leq m(A_i)$, we let

$$O_{i, j} := \left( \bigcup_{l \in \{ 1, \dots, r_i \}; P_{i, l} \in B_{i, j}} (F_{i, l} \cap F^\vee_{i, l}) \right) \cup \left( \bigcup_{l \in \{ 1, \dots, r_{i^\vee} \}; P_{i^\vee, l} \in B_{i^\vee, j}} (F_{i^\vee, l} \cap F^\vee_{i^\vee, l}) \right) \subseteq V,$$

\noindent where we identify $B_{i, j}$ with $B_{i^\vee, j}$ after a suitable numbering of the orbits.\footnote[1]{Since $Q_i = Q_{i^\vee}$, we have $A_{i^\vee} = \alpha_i \cdot A_i$ for some $\alpha_i \in \c^*$. Thus we have $B_{i^\vee, j} = \alpha_i \cdot B_{i, j}$ for a suitable numbering of the orbits. Note, that this correspondence is well-defined, since two such $\alpha_i$ differ by an element of $G(A_i) = G(A_{i^\vee})$, which leave the $B_{i^\vee, j}$ invariant.} In both cases, we let

$$O_0 := V \backslash \left( \bigcup_{i \not\in \mathfrak{E}_D \cup \mathfrak{E}^\vee_{D^\vee}, \ j \in \{ 1, \dots, m(A_i) \} } O_{i, j} \right).$$

\noindent Then the following hold:
\begin{itemize}
\item[(1)] The set $O_0$ is invariant under the action of $\Aut(V)$ and it contains the big orbit $O$. Moreover, the subsets $O_{i, j}$ are invariant under the action of $\Aut(V)$.
\item[(2)] If $F$ is a feather attached to some $C_i$ with $i \in \mathfrak{E}_D \cup \mathfrak{E}^\vee_{D^\vee}$, then any point of $F \backslash D$ is contained in the big orbit $O$.
\item[(3)] If $Q(X, D) \neq Q(X^\vee, D^\vee)$, then for the fixed point set $F(V)$ of the natural action of $\Aut(V)$ on $V$ we have

$$\bigcup_{i \not\in \mathfrak{E}_D \cup \mathfrak{E}^\vee_{D^\vee}, \ d(A_i) = 1, \ j \in \{ 1, \dots, r_i \}} (F_{i, j} \cap F^\vee_{i, j}) \subseteq F(V).$$

\item[(4)] If $r_i > 0$ holds for a unique $i \in \{ 3, \dots, n - 1 \}$, then $O_0$ and the $O_{i, j}$ form the orbit decomposition of the natural action of $\Aut(V)$ on $V$. Moreover, equality in (3) holds.
\end{itemize}
\end{thm}

The most important step for proving Theorem \ref{FixedPointOfSurface} is to show that fibered modifications map the subsets $O_0$ and $O_{i, j}$ onto similar subsets in the image. In other words, we prove the following crucial lemma:

\begin{lemma}\label{FixedPointOfJonquieresAut} Let $V$ be as in $(*)$, $(X, D)$ and $(X', D')$ two standard completions of $V$ and let $\psi: (X, D) \dasharrow (X', D')$ be a fibered modification. Further, let $A_i \subseteq C_i \backslash (C_{i - 1} \cup C_{i + 1})$, $A'_i \subseteq C'_i \backslash (C'_{i - 1} \cup C'_{i + 1})$, $O_{i, j} \subseteq X \backslash D \cong V$ and $O'_{i, j} \subseteq X' \backslash D' \cong V$ be as in Theorem \ref{FixedPointOfSurface}. Then the following hold:
\begin{itemize}
\item[(1)] $\varphi$ maps $O_0$ onto $O'_0$ and $O'_{i, j}$ onto $O'_{i, j}$.
\item[(2)] If $F$ is a feather which is attached to some $C_i$ with $i \in \mathfrak{E}_D$, then $F \backslash D$ is contained in $O$.
\item[(3)] If $Q(X, D) \neq Q(X^\vee, D^\vee)$, then for any $i, j$ with $i \not\in \mathfrak{E}_D \cup \mathfrak{E}^\vee_{D^\vee}$, $d(A_i) = 1$, $j \in \{ 1, \dots, r_i \}$, we have

$$\varphi(F_{i, j} \cap F^\vee_{i, j}) = F'_{i, j} \cap {F'}^\vee_{i, j}.$$

\item[(4)] If $r_i > 0$ holds for a unique $i \in \{ 3, \dots, n - 1 \}$, then $F_{i, j} \backslash (F_{i, j}^\vee \cup D)$ is contained in $O$ if $4 \leq i \leq n - 2$ and, moreover, $F_{i, j} \backslash D$ is contained in $O$ if $i = 3$ or $i = n - 1$.
\end{itemize}
\end{lemma}

\begin{bew} We decompose the map $(X, D) \to (\p^1 \times \p^1, C_0 \cup C_1 \cup C_2)$ into single blowups

$$(X, D) = (X_N, D_N) \stackrel{\pi_N}{\to} \cdots \stackrel{\pi_2}{\to} (X_1, D_1) \stackrel{\pi_1}{\to} (X_0, D_0) = (\p^1 \times \p^1, C_0 \cup C_1 \cup C_2),$$

\noindent according to the algorithm given in \ref{CoordinatesGeneralCase} and we use the coordinate systems $(u_i, v_i)$ and $(x_i, y_i)$ introduced in \ref{CoordinatesGeneralCase}. Similarly, we introduce the corresponding coordinate systems for $(X', D') \to \p^1 \times \p^1$. Then $\psi$ is given on $\p^1 \times \p^1 \backslash (C_0 \cup C_1)$ by

$$\psi(x_0, y_0) = (ax_0 + P(y_0), by_0 + c) \quad \text{for some} \quad a, b \in \c^*, c \in \c, P[y_0] \in \c[y_0].$$

We may assume that the blowup $(X, D) \to (\p^1 \times \p^1, C_0 \cup C_1 \cup C_2)$ is centred in $(0, 0)$. It follows that $P(0) = 0 = c$ and thus

$$\psi(x_0, y_0) = (ax_0 + y_0\tilde{P}(y_0), by_0) \quad \text{for some} \quad \tilde{P}(y_0) \in \c[y_0].$$

We fix an index $j \in \{ 3, \dots, n - 1 \}$. To create the inner components of $D$ and $D'$ we follow the algorithm given in \ref{CoordinatesGeneralCase}. Following step (1) we introduce on $(X_i, D_i)$ coordinate systems $(s_i, t_i)$, $i = 0, \dots, n - 2$, starting with $(s_0, t_0) := (x_0, y_0)$ and such that $(s_{n - 2}, t_{n - 2}) = (u_j, v_j)$. These coordinate systems are related either by $(s_i, t_i) = (s_{i + 1}, t_{i + 1})$ or by $(s_i, t_i) = (s_{i + 1}, s_{i + 1}t_{i + 1})$ or by $(s_i, t_i) = (s_{i + 1}t_{i + 1}, t_{i + 1})$ (and similarly for $(X', D') \to Q$). The first relation arises if we create an inner component at infinity. As in the proof of Theorem \ref{SurfacesRigidExtendedDivisor} we refer to the last two transformations as to transformations of type 1 and type 2 respectively.
 
In the following we denote the lift $(X_i, D_i) \dasharrow (X'_i, D'_i)$ of $\psi$ by $\Psi$, keeping in mind on which intermediate surface $\Psi$ lives.\\ 

\noindent \textbf{Claim 1:} In the coordinates $(u_j, v_j)$ and $(u'_j, v'_j)$, $\Psi$ has the form

\begin{eqnarray}\label{FormOfAutomorphism}
\Psi(u_j, v_j) = ( \alpha u_j(1 + u_j^kv_j^lR(u_j, v_j)), \beta v_j(1 + u_j^kv_j^lS(u_j, v_j)) ),
\end{eqnarray}

\noindent where both components are rational functions, expressed by some power series $R, S \in \c[[u_j, v_j]]$, and $k \geq 0$, $l \geq 1$. Moreover, $k = 0$ holds if and only if we perform no blowup of type 1, that is of type $(s_i, t_i) = (s_{i + 1}, s_{i + 1}t_{i + 1})$ except for the first one in $(0, 0) \in C_2 \backslash C_1$.\\

\noindent \textit{Proof of Claim 1:} We show this by induction on $i$. Since $(s_0, t_0) = (s_1, t_1z_1)$, we have for $i = 1$ 

\begin{eqnarray*}
\Psi(s_1, t_1) &=& \left( as_1 + s_1t_1P(s_1t_1), \frac{bs_1t_1}{as_1 + s_1t_1P(s_1t_1)} \right),\\
&=& \left( as_1(1 + t_1\tilde{P}(s_1, t_1)), \frac{(b/a)t_1}{1 + t_1\tilde{P}(s_1, t_1)} \right)\\ 
&=& ( \alpha s_1(1 + t_1\tilde{P}(s_1, t_1)), \beta t_1(1 + t_1S(s_1, t_1)) ), 
\end{eqnarray*}

\noindent where $\beta = a/b$, $\tilde{P}(s_1, t_1) = a^{-1}P(s_1t_1)$ and the power series $S$ is defined by $1 + t_1S(s_1, t_1) = 1/(1 + a^{-1}P(s_1t_1))$. Thus the claim holds for $i = 1$. In the case $(s_i, t_i) = (s_{i + 1}t_{i + 1}, t_{i + 1})$ we obtain

\begin{eqnarray*}
\Psi(s_{i + 1}, t_{i + 1}) &=& \left( \frac{\alpha}{\beta}s_{i + 1} \cdot \frac{ 1 + s_{i + 1}^kt_{i + 1}^{k + l}R(s_{i + 1}t_{i + 1}, t_{i + 1}) }{ 1 + s_{i + 1}^kt_{i + 1}^{k + l}S(s_{i + 1}t_{i + 1}, t_{i + 1}) }, \beta t_{i + 1}(1 + s_{i + 1}^kt_{i + 1}^{k + l}S(s_{i + 1}t_{i + 1}, t_{i + 1})) \right)\\
&=& ( \tilde{\alpha}s_{i + 1}(1 + s_{i + 1}^kt_{i + 1}^{k + l}\tilde{R}(s_{i + 1}, t_{i + 1}), \beta t_{i + 1}(1 + s_{i + 1}^kt_{i + 1}^{k + l}\tilde{S}(s_{i + 1}, t_{i + 1})) )
\end{eqnarray*}

\noindent with $\tilde{\alpha} = \alpha/\beta$ and certain power series $\tilde{R}, \tilde{S}$. Similarly, in the case $(s_i, t_i) = (s_{i + 1}, s_{i + 1}t_{i + 1})$ we obtain

\begin{eqnarray*}
\Psi(s_{i + 1}, t_{i + 1}) &=& ( \alpha s_{i + 1}(1 + s_{i + 1}^{k + l}t_{i + 1}^l\tilde{R}(s_{i + 1}, t_{i + 1}), \tilde{\beta} t_{i + 1}(1 + s_{i + 1}^{k + l}t_{i + 1}^l\tilde{S}(s_{i + 1}, t_{i + 1})) )
\end{eqnarray*}

\noindent with $\tilde{\beta} = \beta/\alpha$ and certain power series $\tilde{R}, \tilde{S}$. Now, looking at the induction step it is obvious that $k = 0$ holds if and only if we perform no blowup of type 1 except for the first one in $(0, 0) \in C_2 \backslash C_1$.\\

Now we start to create feathers on the components $C_3, \dots, C_{j - 1}$. This corresponds to step (3) in \ref{CoordinatesGeneralCase}. It is easy to verify the following claim:\\

\noindent \textbf{Claim 2:} Let $m < j$. Replacing $(u_j, v_j)$ by $(u'_j, v'_j) = (T_{mj} \circ T_k^{\text{Bl.up}} \circ T_{mj}^{-1})(u_j, v_j)$, the map $\Psi$ has the same form as is (\ref{FormOfAutomorphism}) with the same values for $k$ and $l$ (but with possibly other power series $R$ and $S$).\\
 



In particular, Claim 2 implies the following: given any $j \in \{ 3, \dots, n - 1 \}$, then $\Psi: X_k \dasharrow X'_k$ has the form (\ref{FormOfAutomorphism}) on any intermediate surface $X_k$ where only feathers with mother components $C_\tau$, $\tau < i$, are created.

Now we want to observe the action of $\Psi$ on the feathers attached to $C_j$. To create the feathers on the component $C_j$ we blow up $r_j$ points on $C_j \backslash (C_{j - 1} \cup C_{j + 1})$ (and similarly for $C'_j$). In the coordinates $(u_j, v_j)$ and $(u'_j, v'_j)$, respectively, the base point sets $A_j$ and $A'j$ can be written as $A_j = \{ P_{j, i} = (\gamma_i, 0) \mid 1 \leq i \leq r_j \}$ and $A'_j = \{ P'_{j, i} = (\gamma_i, 0) \mid 1 \leq i \leq r_j \}$, respectively. Since $\psi$ lifts to $X \dasharrow X'$ we have $\Psi(A_j) = A'_j$. The restriction

$$\Psi\vert_{C_j \backslash (C_{j - 1} \cup C_{j + 1})}: C_j \backslash (C_{j - 1} \cup C_{j + 1}) \to C'_j \backslash (C'_{j - 1} \cup C'_{j + 1})$$

\noindent induces the multiplication $x \mapsto \alpha \cdot x$, \ie $\alpha \in G(A_j)$. We consider points $P_{j, s}$ and $P'_{j, t}$ such that $\Psi(P_{j, s}) = P'_{j, t}$, \ie such that $\alpha \cdot \gamma_s = \gamma_t$ holds. After the blowup in $P_{j, s}$ and $P'_{j, t} = \Psi(P_{j, s})$, respectively, we introduce affine coordinates $(u, v)$ and $(u', v')$, respectively, via $(x_j, y_j) = (uv + \gamma_s, v)$ and $(x_j, y_j) = (u'v' + \gamma_t, v')$, respectively. We represent the arguments of $\Psi$ in the coordinates $(u, v)$ and the images in the coordinates $(u', v')$. In these coordinates we have $F_s \backslash D = \{ v = 0 \}$ and $F'_t \backslash D' = \{ v' = 0 \}$, and the lift of $\Psi$ has the form

\begin{eqnarray}\label{Lifting}
\Psi(u, v) &=& \left( \frac{\alpha}{\beta}u + v^{l - 1}\tilde{R}(u, v), \beta v \cdot (1 + v^l\tilde{S}(u, v)) \right)
\end{eqnarray}
 
\noindent for certain power series $\tilde{R}, \tilde{S}$. Thus we see that $\tilde{\Psi} \vert_{F_s \backslash D}: F_s \backslash D \to F'_t \backslash D'$ is given by

$$\Psi(u, 0) = \begin{cases} \left( \frac{\alpha}{\beta}u + \tilde{R}(u, 0), 0 \right) &, l = 1 \\ \left( \frac{\alpha}{\beta}u, 0 \right) &, l \geq 2.\end{cases}$$ 

\noindent Since $(u, v) = (0, 0)$ and $(u', v') = (0, 0)$, respectively, are precisely the intersection points $F_{j, s} \cap F^\vee_{j, s}$ and $F'_{j, t} \cap {F'}^\vee_{j, t}$, respectively, we obtain that $\Psi(F_s \cap F_s^\vee) = F'_t \cap {F'}_t^\vee$ holds if $l \geq 2$. The next claim finishes the proof of assertion (1):\\

\noindent \textbf{Claim 3:} $l \geq 2$ holds if and only if $C_j$ is not an exceptional component of $D$.\\

\noindent \textit{Proof of Claim 3:} We consider the construction of the extended divisor $D_{\text{ext}}$ starting from the zigzag $C_0 \triangleright C_1 \triangleright C_2$ on $Q$. We denote the exponents $k, l$ in the $i$-th lift of $\Psi$, that is, the exponents in the expression of $\Psi(s_i, t_i)$, now by $k_i, l_i$, $i = 0, \dots, n - 2$. The first blowup in $(0, 0) \in C_2 \backslash C_1$ is of type 1 and we obtain $k_1 = 0$ and $l_1 = 1$. During further blowups we can first perform as many blowups as possible of type 2. After this, we get a certain intermediate boundary divisor $D_p = C_0 \triangleright C_1 \triangleright C_2 \triangleright C_{\mu_3} \triangleright \cdots \triangleright C_{\mu_{p + 1}} \triangleright C_n$, $2 \leq p \leq n - 2$, with dual the graph $\Gamma_{D'} = [[0, 0, -2, \dots, -2, -1, -p ]]$. Note, that the components $C_{\mu_3}, \dots, C_{\mu_{p + 1}}$ are all exceptional components by definition. During these steps we obtain $k_i = 0$, $l_i = 1$ (see induction step). Here the right-hand boundary component gives the $t_i$-axis and the next-to-last one gives the $s_i$-axis. If the feathers $F_{j, k}$ are attached to $C_{\mu_{p + 1}}$, we need no more transformation of type 2 and we end up in $l = l_{n - 2} = 1$. But if $j \not\in \mathfrak{E}_D$, we need at least one more transformation of type 1 and thereafter at least one transformation of type 2 (since the base points of the feathers $F_{j, k}$ lie on the $s_{n - 2}$-axis). Since transformations of type 1 give $k_{i + 1} = k_i + l_i$ and $l_{i + 1} = l_i$ and transformations of type 2 give $k_{i + 1} = k_i$ and $l_{i + 1} = k_i + l_i$, it follows that after performing all blowups we have $k_i \geq 1$ and $l_i \geq 2$. In particular we obtain $l_{n - 2} \geq l_i \geq 2$.\\

For (2), let $F$ be a feather that is attached to an exceptional component $C_i$, $i \in \mathfrak{E}_D$, of $D$. We consider once again the blowup process $X_{n - 2} \to X_0 = Q$, which produces the boundary divisor $D$ and stop it at a suitable intermediate surface. Since $F$ is attached to an exceptional component, we have to stop at a step where the intermediate boundary divisor $\tilde{D}$ has the dual graph $\Gamma_{\tilde{D}} = [[0, 0, -2, \dots, -2, -1, -p ]]$, such that $C_i$ is the proper transform of the unique $(-1)$-curve $\tilde{C}_{p + 1}$ of $\tilde{D}$ (and similarly we proceed for $X' \to Q$). We consider the lift of the fibered map $\psi_a(x_0, y_0) = (x_0 + ay_0, y_0)$, $a \in \c$, on $(\tilde{X}, \tilde{D})$. As showed above, such maps lift to fibered maps $(X, D) \dasharrow (X', D')$ (and in particular to any intermediate surface). We consider again the coordinates $(s_i, t_i)$, introduced in \ref{CoordinatesGeneralCase} (1). On the surface $(\tilde{X}, \tilde{D})$ we have the coordinates $(s_p, t_p)$ with $\tilde{C}_{p + 2} \backslash \{ \infty \} = \{ s_p = 0 \}$, $\tilde{C}_{p + 1} \backslash \tilde{C}_p = \{ t_p = 0 \}$ and $(x_0, y_0) = (s_0, t_0) = (s_pt_p^{p - 1}, s_pt_p^p)$, and the lift $\Psi_a: (\tilde{X}, \tilde{D}) \dasharrow (\tilde{X}', \tilde{D}')$ of $\psi_a$ is given by

$$\Psi_a(s_p, t_p) = \left( s_p(1 + at_p)^p, \frac{t_p}{1 + at_p} \right).$$

\noindent If $F$ has the base point $(s_p, t_p) = (\gamma, 0)$, we introduce the coordinates $(u, v)$ via $(s_p, t_p) = (uv + \gamma, v)$. Similarly, since $\Psi_a(\gamma, 0) = (\gamma, 0)$, we introduce on $(\tilde{X}', \tilde{D}')$ the coordinates $(u', v')$. These coordinates satisfy $F \backslash D = \{ u = 0 \}$ and $F' \backslash D' = \{ u' = 0 \}$, and the restriction $F \backslash D \to F' \backslash D'$ of the lift of $\psi_a$ is given by

$$\Psi_a(u, 0) = (u + pa\gamma, 0).$$

\noindent Since the coordinates $(u', v')$ play the same role for $(X', D')$ as $(u, v)$ play for $(X, D)$, it follows that $F \backslash D \subseteq O$.

Assertion (3) is an immediate consequence of (1) since $d(A_i) = 1$ and $Q(X, D) \neq Q(X^\vee, D^\vee)$ imply that the sets $O_{i, j}$ are reduced to single points.

Finally, we assume that $r_i > 0$ and $r_j = 0$ for $j \neq i$. First, it is not hard to see that in this case, if the feathers are attached to an exceptional component, then they are attached to $C_{n - 1}$. Similarly as in (\ref{ArbitraryMultiples}), the torus element $t = (a, b) \in \t$ gives the map $(u_i, v_i) \mapsto (a^{p_i}b^{q_i}u_i, a^{-p_{i + 1}}b^{-q_{i + 1}}v_i)$, since $v_i = u_{i + 1}^{-1}$ and no feather is created by a blowup (note, that $t$ defines an isomorphism $t: (X, D) \stackrel{\sim}{\to} (X, D)$). Assume now that the base points $P_{i, s}$ and $P_{i, t}$ are contained in the same $G(A_i)$-orbit. Then there exists an appropriate $t$ with $t(P_{i, s}) = P_{i, t}$. Blowing up in $A_i$, $t$ induces the map

$$t: F_{i, s} \backslash D \to F_{t, s} \backslash D, \quad u \mapsto a^{p_i + p_{i + 1}}b^{q_i + q_{i + 1}}u.$$
  
\noindent Since 

$$\left| \begin{array}{cc} p_i & p_i + p_{i + 1} \\ q_i & q_i + q_{i + 1} \end{array} \right| = \left| \begin{array}{cc} p_i & p_{i + 1} \\ q_i & q_{i + 1} \end{array} \right| \neq 0,$$ 
  
\noindent for any given $\alpha \in \c^*$, the elements $a, b$ can be chosen, in addition, in such a way that $a^{p_i + p_{i + 1}}b^{q_i + q_{i + 1}} = \alpha$. This shows that the orbit of any point in $F_{i, j} \backslash (F_{i, j}^\vee \cup D)$ is infinite and is thus contained in $O$. 

Consider now the case $i = n - 1$, \ie all feathers are attached to $C_{n - 1}$. Let $F$ be such a feather. By Claim 3, the fibered maps $\psi_a(x_0, y_0) = (x_0 + ay_0, y_0)$, $a \in \c$, induce non-trivial translations on $F \backslash D$ for $a \neq 0$. Thus $F \backslash D$ is contained in $O$. Similarly, passing to $(X^\vee, D^\vee)$ it follows that $F \backslash D \subseteq O$, if $F$ is attached to $C_3$. Further, equality in (3) holds due to the fact that if $r_i > 0$ holds for a unique $i$, we have $Q(X, D) = Q(X^\vee, D^\vee)$ if and only if $D$ is of the type $[[0, -1, -2, -1 - r_3, -2]]$. However, in this case $\Aut(V)$ acts transitively on $V$, since $C_3$ is an exceptional component of $D$. Hence (4) follows.    
\end{bew}

Another important step in the proof of Theorem \ref{FixedPointOfSurface} is to show that a reversion $\varphi$ maps the subsets $O_{i, j} \subset V$ onto similar subsets in the image of $\varphi$.

\begin{lemma}\label{ReversionFixIntersectionPoints} Let $V$ be a smooth Gizatullin surface as in Theorem \ref{FixedPointOfSurface} and let $(X, D)$ be a standard completion of $V$. Further, let $\varphi: (X, D) \dasharrow (X', D')$ be a reversion and let $O_{i, j}$ and $O'_{i, j}$ denote the corresponding subsets of $X$ and $X'$ defined in Theorem \ref{FixedPointOfSurface}.
\begin{itemize}
\item[(a)] If $Q(X, D) \neq Q(X^\vee, D^\vee)$, then $\varphi$ maps $O_{i, j}$ onto $O_{i^\vee, j}'$, where we identify $A_i$ with $A'_{i^\vee, j}$ via the Matching Principle.
\item[(b)] If $Q(X, D) = Q(X^\vee, D^\vee)$, then $\varphi$ maps $O_{i, j}$ onto $O_{i, j}'$ for all $i \not\in \mathfrak{E}_D \cup \mathfrak{E}^\vee_{D^\vee}$, $i < \lfloor \frac{n}{2} \rfloor + 1$.
\end{itemize}
\end{lemma}

\begin{bew} 

By considering a resolution $(X, D) \leftarrow (Z, B) \to (X', D')$ of $\varphi$ it is clear that the proper transforms of $F_{i, j}$ and $F_{i, j}^\vee$, respectively, become ${F^\vee}'_{i^\vee, j}$ and $F'_{i^\vee, j}$, respectively, under a suitable numbering of the feathers. This yields, by construction of the $O_{i, j}$, that $\varphi$ maps $O_{i, j}$ onto $O_{i^\vee, j}$ if $Q(X, D) \neq Q(X^\vee, D^\vee)$, and $O_{i, j}$ onto $O_{i, j}$ if $Q(X, D) = Q(X^\vee, D^\vee)$.
\end{bew}

Let us now prove Theorem \ref{FixedPointOfSurface}.

\bigskip
\noindent \textit{Proof of Theorem \ref{FixedPointOfSurface}:} Let $\varphi \in \Aut(V)$. We extend $\varphi$ to a birational map $\varphi: (X, D) \dasharrow (X, D)$, which belongs either to $\Aut(X, D)$ or admits a decomposition $\varphi = \varphi_m \circ \cdots \circ \varphi_1$, where each $\varphi_i$ is either a reversion or a fibered modification. Thus $\varphi$ decomposes into a sequence

$$(X, D) = (X_0, D_0) \stackrel{\varphi_1}{\dasharrow} (X_1, D_1) \stackrel{\varphi_2}{\dasharrow} \cdots \stackrel{\varphi_m}{\dasharrow} (X_m, D_m) = (X, D).$$

Let now $O^{(k)}_{i, j}$ be the corresponding subsets of $V = X_k \backslash D_k$ as defined in Theorem \ref{FixedPointOfSurface}. In the following we have to distinguish two cases:

If $\varphi_k: (X_{k - 1}, D_{k - 1}) \dasharrow (X_k, D_k)$ is a fibered modification, then by Lemma \ref{FixedPointOfJonquieresAut}, then $\varphi_k$ maps the subsets $O^{(k - 1)}_{i, j}$ onto $O^{k}_{i, j}$. 

If $\varphi_k: (X_{k - 1}, D_{k - 1}) \dasharrow (X_k, D_k)$ is a reversion, then by Lemma \ref{ReversionFixIntersectionPoints} the subsets $O^{(k - 1)}_{i, j}$ are mapped onto $O^{(k)}_{i^\vee, j}$, if $Q(X, D) \neq Q(X^\vee, D^\vee)$, and onto $O^{(k)}_{i, j}$, if $Q(X, D) = Q(X^\vee, D^\vee)$.

Now, in the case $Q(X, D) \neq Q(X^\vee, D^\vee)$, the number of reversions occuring in the decomposition of $\varphi$ is even. Hence, in both cases $\varphi(O_{i, j}) = O_{i, j}$. This gives assertion (1).

Let $F$ be a feather which is attached to $C_i$ with $i \in \mathfrak{E}_D$. Then by Lemma \ref{FixedPointOfJonquieresAut} (2), $F \backslash D$ is contained in $O$. Reversing $(X, D)$ yields that $F \backslash D$ is contained in $O$, if $F$ is attached to some $C_i$ with $i \in \mathfrak{E}^\vee_{D^\vee}$. This gives (2).

Assertion (3) follows immediately from Lemma \ref{FixedPointOfJonquieresAut} (3) and assertion (4) follows from Lemma \ref{FixedPointOfJonquieresAut} (4). This ends the proof.

\begin{picture}(2,2)
\put(433,0) {\line(1,0){8}}
\put(441,0) {\line(0,1){8}}
\put(441,8) {\line(-1,0){8}}
\put(433,8) {\line(0,-1){8}}
\end{picture}
\\

\begin{ex} Consider a smooth Gizatullin surface $V$ with a standard completion

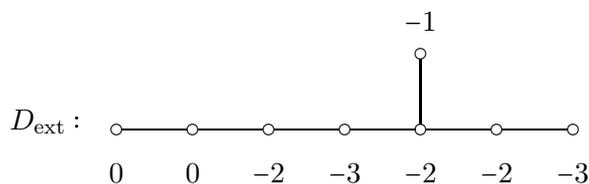
\vspace{15pt}
$$D_{\text{ext}}: \quad \cu{0} \lin \cu{0} \lin \cu{-2} \lin \cu{-3} \lin \cu{-2} \nlin \cshiftup{-1}{} \lin \cu{-2} \lin \cu{-3}$$

The exceptional components are $C_3$ and $C_5$ and thus $\mathfrak{E}_D = \{ 3, 5 \}$. Indeed, the boundary divisor arises in the following way:

$$\cou{C_0}{0} \lin \cou{C_1}{0} \lin \cou{C_2}{0} \quad \leftarrow \quad \cou{C_0}{0} \lin \cou{C_1}{0} \lin \cou{C_2}{-2} \lin \cou{C_3}{-2} \lin \cou{C_5}{-1} \lin \cou{C_6}{-3} \quad \leftarrow \quad \cou{C_0}{0} \lin \cou{C_1}{0} \lin \cou{C_2}{-2} \lin \cou{C_3}{-3} \lin \cou{C_4}{-1} \lin \cou{C_5}{-2} \lin \cou{C_6}{-3}$$

Further, a direct computation yields that $\mathfrak{E}_{D^\vee} = \{ 5 \}$. Hence, by Theorem \ref{FixedPointOfSurface} there is a unique fixed point $p := F \cap F^\vee$, where $F$ is the unique feather of $D_{\text{ext}}$. Moreover, the action of $\Aut(V)$ has exactly two orbits, namely $O = V \backslash \{ p \}$ and the fixed point $\{ p \}$.
\end{ex}

\begin{rem} The proof of Theorem \ref{FixedPointOfSurface} only shows that the subsets $O_{i, j}$ are $\Aut(V)$-invariant, but it does not say anything about the finite subset $\bigcup_{i \neq j} F_{i, \rho} \cap F^\vee_{j, \sigma}$. The reason is that we have no description of these intersection points in the coordinates $(x_i, y_i)$ introduced in \ref{CoordinatesGeneralCase}. Theoretically, the algorithm given in \ref{CoordinatesGeneralCase} allows us to compute explicit presentations of all dual feathers $F_{j, \sigma}$ and thus to compute these intersection points. However, this would mean an immense computation and so we omit it.  
\end{rem}

\begin{rem} If $V$ admits more than two families of feathers, then the subsets $O_{i, j}$ are not orbits of $\Aut(V)$ in general and there may be a strict inclusion in Theorem \ref{FixedPointOfSurface} (2). The problem is obvious: we have only two parameters, which induce motions $x \mapsto \alpha_i \cdot x$ on the $*$-components $C_i$, $4 \leq i \leq n - 2$ (these parameters are given by $a, b \in \c^*$ in the fibered modifications $\psi(x_0, y_0) = (ax_0 + y_0Q(y_0), by_0)$). Moreover, it is difficult to control these motions since they cannot be read out directly from the form of the extended divisor.

\end{rem}

\begin{rem} Let $V$ be a normal affine variety. We denote by $\SAut(V)$ the subgroup of $\Aut(V)$ generated by all algebraic subgroups isomorphic to the additive group $\mathbb G_a$. Then a point $x \in V$ is called \textit{flexible}, if the tangent space $T_xV$ is spanned by the tangent vectors of the orbits $H.x$ of one-parameter unipotent subgroups $H \subseteq \Aut(V)$. Moreover, we call $V$ \emph{flexible} if every point $x \in V_{\text{reg}}$ is flexible.
 
If $V$ is a flexible affine variety then $\SAut(V)$ acts transitively on $V_{\text{reg}}$. Moreover, in \cite{AFKKZ} it was shown that for a normal affine variety $V$ of dimension $\geq 2$ the following conditions are equivalent (see \cite{AFKKZ}, Theorem 0.1):

\begin{itemize}
\item[(1)] $X$ is flexible.
\item[(2)] $\SAut(V)$ acts transitively on $V_{\text{reg}}$.
\item[(3)] $\SAut(V)$ acts infinitely transitively on $V_{\text{reg}}$.
\end{itemize}

The last condition means that for any collection of points $\{ P_1, \dots, P_k \}$ and $\{ Q_1, \dots, Q_k \}$ in $V_{\text{reg}}$ there exists an automorphism $\varphi \in \SAut(V)$ such that $\varphi(P_i) = Q_i$.\\
Therefore, instead of "flexible" we sometimes say "infinitely transitive". Similarly, $V$ is called \textit{stably infinitely transitive} if $V \times \a^m$ is infinitely transitive for some $m \geq 0$.

Theorem \ref{FixedPointOfSurface} shows that smooth Gizatullin surfaces $V$ with a distinguished and rigid extended divisor are \textit{not flexible} in general. But the following question arises:\\

\noindent \textbf{Question/Problem:} Let $V$ be a smooth Gizatullin surface with a distinguished and rigid extended divisor. Is $V$ stably infinitely transitive?
\end{rem}

\begin{rem}\label{OrevkovCurves}
Recently, M. Kh. Gizatullin informed me that V. I. Danilov discovered 1973 a class of curves, today known as Orevkov curves, which give counterexamples to Gizatullins conjecture (see the review to \cite{MS}). Unfortunately, these counterexamples were never published.

Orevkov curves are plane curves $F_n \subseteq \p^2$, $n \geq 1$, of degree $f_{2n - 1}$, where $f_k$ is the $k$-th Fibonacci number (starting with $f_0 = 0$). These curves admit a unique singular point. We consider the affine surfaces $V_n := \p^2 \backslash F_n$. Resolving the singularity yields that $V_n$ is completable by a zigzag. Denoting by $l_n$ the lenght of the zigzag of $V_n$ in standard form, it can be shown that the extended divisor is distinguished and rigid and moreover, that it admits a feather which is attached to an inner component $C_{i_n}$, which is not exceptional, if $n \geq 4$ (\ie, if the degree of $F_n$ is at least $13$). By Theorem \ref{FixedPointOfSurface}, the action of $\Aut(V_n)$ admits a fixed point.
\end{rem}

\subsection{The amalgamated product structure of the automorphism group}\label{AmalgamatedProductStructure}

The automorphism groups of surfaces considered in Theorem \ref{SurfacesRigidExtendedDivisor} can be represented as amalgamated products of certain subgroups.

\begin{cor}\label{StructureOfAutomorphismGroup} Let $V$ be as in Theorem \ref{SurfacesRigidExtendedDivisor} and let $(X, D)$ be a $1$-standard completion of $V$. We choose a fixed $\a^1$-fibration $\pi: V \to \a^1$ and consider the corresponding $\a^1$-fibration $\pi^\vee: V \to \a^1$, induced by the reversion $\psi: (X, D) \dasharrow (X^\vee, D^\vee)$ with center $p \in C_0 \backslash C_1$.  
\begin{itemize}
\item[(1)] If $\F_V$ consists of a vertex and a loop, \ie $(X, D) \cong (X^\vee, D^\vee)$, then the automorphism group of $V$ is the free product of $A = \langle \Aut(X, D), \psi \rangle$ and $J = \Aut(V, \pi)$, amalgamated over their intersection $A \cap J = \Aut(X, D)$:

$$\Aut(V) = A \star_{A \cap J} J.$$

\item[(2)] Let $\F_V$ be of the form $\begin{xy} \xymatrix{ [(X, D)] \ \bullet \ar@{<->}[r] & \bullet \ [(X^\vee, D^\vee)] } \end{xy}$. We denote by $A$ the subgroup corresponding to the edge and by $J$ and $J^\vee$ the subgroups $J := \Aut(V, \pi)$ and $J^\vee := \Aut(V, \pi^\vee)$. Identifying $J \supseteq A \subseteq J^\vee$ we have

$$\Aut(V) = J \star_A J^\vee.$$

\end{itemize}
\end{cor}

\begin{bew} Let $\F_V$ have a loop. First we show that $\Aut(V)$ is generated by $\Aut(X, D), \psi$ and $J$. We can extend every automorphism $\varphi \in \Aut(V)$ to a birational map $\varphi: (X, D) \dasharrow (X, D)$. Then $\varphi$ either belongs to $\Aut(X, D)$ or it can be decomposed as

$$\varphi = \varphi_n \circ \cdots \circ \varphi_1: (X, D) = (X_0, D_0) \stackrel{\varphi_1}{\dasharrow} (X_1, D_1) \stackrel{\varphi_2}{\dasharrow} \cdots \stackrel{\varphi_m}{\dasharrow} (X_m, D_m) = (X, D),$$

\noindent where every $\varphi_i$ is a fibered modification or a reversion. Since every $1$-standard completion of $V$ is isomorphic to $(X, D)$, we can assume that $(X_i, D_i) = (X, D)$. Therefore, every element $\varphi_i$ can be considered as an element of $\Aut(V)$. If $\varphi_i$ is a fibered modification, then we have $\varphi_i \in J$. But if $g_i$ is a reversion, then it can be written as $\varphi_i = \alpha_i \psi \beta_i$ with certain $\alpha_i, \beta_i \in \Aut(X, D)$, since $\Aut(X, D)$ acts transitively on $C_0 \backslash C_1$. Thus we obtain $\Aut(V) = \langle \Aut(X, D), \psi, J \rangle = \langle A, J \rangle$.

Now we write any $\varphi \in \Aut(V)$ as $\varphi = a_n \circ j_n \circ \cdots \circ a_1 \circ j_1$ with $a_i \in A \backslash J$ and $j_i \in J \backslash A$. Then $a_i$ is a product of reversions which is not an isomorphism, and $j_i$ is a fibered modification. Theorem \ref{FactorizationOfBirationalMaps} then yields that $\varphi \not\in \Aut(X, D)$. Thus it follows that $\Aut(V) = A \star_{A \cap J} J$.

Assertion (2) follows immediately from Remark \ref{ExTreeWithTwoVertices}.
\end{bew}

\section{Remarks on the singular case}

Some of our results can be generalized to the case of singular Gizatullin surfaces. First, we have to generalize the notion of a $*$-component.

\begin{df}\label{StarComponentGeneral}
\begin{itemize}
\item[(1)] For a general feather $F$ with dual graph

$$\Gamma_F: \quad \cu{B} \lin \cu{D_1} \lin \dots \lin \cu{D_k}$$

\noindent and bridge curve $B$ we call $D_k$ the \emph{tip component of} $F$
\item[(2)] The component $C_i$ is called a \emph{$*$-component} if 
\begin{itemize}
\item[(i)] $D^{\geq i + 1}_{\text{ext}}$ is not contractible and
\item[(ii)] $D^{\geq i + 1}_{\text{ext}} - F_{j, k}$ is not contractible for every feather $F_{j, k}$ of $D^{\geq i + 1}_{\text{ext}}$ such that the tip component of $F_{j, k}$ has mother component $C_\tau$ with $\tau < i$.
\end{itemize}
Otherwise $C_i$ is called a $+$-component.
\end{itemize}
\end{df}

It is easy to see that also in the non-smooth case, that is where feathers can have more than one component, $*$-components appear as a result of an inner blowup of the previous zigzag, while an outer blowup of a zigzag creates a $+$-component. Similarly as in the smooth case we have the following lemma:

\begin{lemma} Let $D_{\text{ext}}$ be the extended divisor of the minimal resolution of singularities of a $1$-standard completion of a certain Gizatullin surface $V$. Suppose that every $C_i$, $3 \leq i \leq n - 1$, is a $*$-component and that there is no feather attached to the component $C_n$. Then every feather $F_{i, j}$ is an $A_k$-feather, that is, every $F_{i, j}$ is contractible and therefore has the dual graph

$$\Gamma_{F_{i, j}}: \quad \cou{-1}{B} \lin \cou{-2}{D_1} \lin \dots \lin \cou{-2}{D_k} \quad ,$$

with $k$ depending on $i$ and $j$. 
\end{lemma}

Note that especially for $A_k$-feathers the mother components of all curves $D_1, \dots, D_k$ coincide since any $A_k$-feather is born by successive blowups of a point on the boundary component it is attached to.

\begin{rem} If every component $C_i$, $3 \leq i \leq n - 1$, of $D_{\text{ext}}$ is a $*$-component and if there are no feathers attached to $C_2$ and $C_n$, then it is easy to see that the same property holds for the extended divisor $D^\vee_{\text{ext}}$ after reversion (with an arbitrary center).
\end{rem}

In the following we will assume that the following condition holds:

\begin{eqnarray*}
(*) && V \ \text{admits a} \ 1-\text{standard completion} \ (X, D) \ \text{such that} \ C_3, \dots, C_{n - 1} \ \text{are} \ *\text{-components} \\
&& \text{and there is no feather attached to} \ C_2 \ \text{and to} \ C_n.
\end{eqnarray*}

\bigskip
We can generalize Theorem \ref{SurfacesRigidExtendedDivisor} to singular Gizatullin surfaces:

\begin{thm}\label{SurfacesRigidExtendedDivisorGeneral} Let $V$ be a Gizatullin surface as in $(*)$ with $n \geq 4$ and let $(X, D)$ be the minimal resolution of singularities of a $1$-standard completion of $V$ . For every $s \geq 0$, we let $A_{i, s} = \{ P_{i, s, 1}, \dots, P_{i, s, r_{i, s}} \} \subseteq C_i \backslash (C_{i - 1} \cup C_{i + 1}) \cong \c^*$, $3 \leq i \leq n - 1$, be the base point set of the feathers $F_{i, j}$, $j \in \{ 1, \dots, r_i \}$, which are $A_s$-feathers. Then the dual graph of $D_{\text{ext}}$ has the form

\bigskip
$$D_{\text{ext}}: \quad \cu{0} \lin \cu{-1} \lin \cu{w_2} \lin \cu{w_3} \nlin \xbshiftup{ \{ F_{3, j} \} }{} \lllin \cu{w_4} \nlin \xbshiftup{ \{ F_{4, j} \} }{} \lllin \dots \lllin \cu{w_{n - 1}} \nlin \xbshiftup{ \{ F_{n - 1, j} \} }{} \lin \cu{w_n}$$

\noindent with $j \in \{ 1, \dots, r_i \}$. Moreover, the following hold:
\begin{itemize}
\item[(1)] For any two $1$-standard completions $(X', D')$ and $(X'', D'')$ of $V$ with $(X' \backslash D', \pi') \cong (X'' \backslash D'', \pi'')$ we have $(X', D', \bar{\pi}') \cong (X'', D'', \bar{\pi}'')$. 
\item[(2)] The graph $\F_V$ has one of the following two forms:

$$\F_V: \begin{xy}
  \xymatrix{
  [(X, D)] \ \bullet \ar@{<->}[r] & \bullet \ [(X^\vee, D^\vee)] \\
  }
\end{xy} \quad \text{or} \quad \F_V: [(X, D)] \ \bullet \rcirclearrowleft .$$

\noindent If $\F_V$ is of the form $\bullet \rcirclearrowleft$, then $D^{\geq 2}$ is a palindrome and there exist elements $\gamma_i \in \c^*$, $3 \leq i \leq n - 1$, such that

$$A_{i^\vee, s} = \gamma_i \cdot A_{i, s} \quad \text{holds for all} \ 3 \leq i \leq n - 1, s \geq 0.$$

\item[(3)] If $r_i > 0$ holds for at most two indices $i \in \{ 3, \dots, n - 1 \}$, then $\F_V$ has the form $\bullet \rcirclearrowleft$ if and only if $D^{\geq 2}$ is a palindrome and there exist elements $\gamma_i \in \c^*$ such that

$$A_{i^\vee, s} = \gamma_i \cdot A_{i, s} \quad \text{holds for} \ 3 \leq i \leq n - 1, s \geq 0.$$

\item[(4)] $\Aut(V)$ is generated by automorphisms of $\a^1$-fibrations if and only if $\F_V$ has no loops except for the case $\Gamma_D = [[0,-1, -2, -2, -2]]$. If $\Gamma_D = [[0,-1, -2, -2, -2]]$, then $V$ is smooth.
\end{itemize}
\end{thm}

\begin{bew} Let $(X', D')$ and $(X'', D'')$ be two $1$-standard completions of $V$ such that $(X' \backslash D', \pi') \cong (X'' \backslash D'', \pi'')$. By Lemma \ref{JonquieresAutomorphisms} such an isomorphism is given by a Jonquieres automorphism of the form $\psi(x_0, y_0) = (ax_0 + y_0P(y_0), by_0)$, $P(y_0) \in \c[y_0]$ of $\a^2 = \f_1 \backslash (C_0 \cup C_1)$. The same computation as in the proof of Theorem \ref{SurfacesRigidExtendedDivisor} shows that $\psi$ lifts to an automorphism of a $1$-standard completion $(X, D)$ of $V$ as well as to the minimal resolution of singularities $(Y, B)$ of $(X, D)$ (notice that $B \cong D$ since $D$ is contained in the regular locus of $X$). Now, a similar computation as in the proof of Theorem \ref{SurfacesRigidExtendedDivisor} shows that $\psi$ describes on $D^{\geq 2}$ (before creating the feathers by blowups of points on the boundary components) the same map as the isomorphism $\tilde{\psi}(x_0, y_0) = (ax_0 + cy_0, by_0)$ with $c = P(0)$. By Lemma \ref{JonquieresAutomorphisms2} we obtain that $(X', D') \cong (X'', D'')$. This shows (1).

The proof of assertion (2), (3) and (4) is the same as in the smooth case. 
\end{bew}

\begin{rem} Corollary \ref{StructureOfAutomorphismGroup} holds as well for arbitrary Gizatullin surfaces satisfying $(*)$, since the proof uses only the structure of the graph $\F_V$.
\end{rem}

\begin{rem} (Toric Surfaces) For a primitive $d$-th root of unity $\zeta$ and integers $d, e \geq 0$ with $0 \leq e < d$ and $\gcd(e, d) = 1$, the cyclic group $\z_d \cong W_d= \langle \zeta \rangle$ acts on $\a^2$ via $\zeta.(x, y) = (\zeta x, \zeta^e y)$. The quotient $V_{d, e} := \a^2//\z_d$ is a normal affine toric surface. Conversely, any normal toric surface non-isomorphic to $\c^* \times \c^*$ or to $\c^* \times \a^1$ arises in this way (see \cite{FKZ3}, 1.8). The surfaces $V_{d, e}$ are Gizatullin surfaces and we describe their extended divisors as well as their automorrphism group. Let $(X, D)$ be a standard completion. By \cite{FKZ2} Lemma 2.20, $D_{\text{ext}}$ is linear. The boundary divisor $D$, up to reversion, has the following dual graph

$$\Gamma_D: \quad \co{0} \lin \co{0} \lin \boxo{\frac{d - e}{d}} \quad ,$$

\noindent where a box with weight $\frac{a}{b}$ abbreviates the zigzag

$$\co{-k_1} \lin \dots \lin \co{-k_n} \quad ,$$

\noindent such that $[k_1, \dots, k_n] = b/a$ (see \cite{FKZ3}, Lemma 3.12). Since $D_{\text{ext}}$ admits a unique feather $F$ and since it can be contracted to $[[0, 0, 0]]$, it is not hard to see that the feather has the dual graph

$$\Gamma_F: \quad \co{-1} \lin \boxo{\frac{e}{d}} \quad .$$

\noindent Reversing $(X, D)$ we obtain a completion $(X^\vee, D^\vee)$ with dual graphs

$$\Gamma_D: \quad \co{0} \lin \co{0} \lin \boxo{\frac{d - e'}{d}} \quad \text{and} \quad \Gamma_{F^\vee}: \quad \co{-1} \lin \boxo{\frac{e'}{d}} \quad ,$$

\noindent where $e'$ is the unique integer with $0 \leq e' < d$ and $ee' \equiv 1 \mod d$. Now, passing to $1$-standard completions we obtain that there are also at most two possible dual graphs of an extended divisor of $V$. Let us now fix a $1$-standard completion $(X, D)$ of $V$. Since the bridge curve of the feather is the only component of $D_{\text{ext}}$ with self-intersection number $-1$, all blowups in $X \to \f_1$ are inner, except for the first one (which creates the tip component of the feather). Using this, the same comutation as in the proof of Theorem \ref{SurfacesRigidExtendedDivisor} yields that $\Aut(X, D)$ acts transitively on $C_0 \backslash C_1$. Hence, concluding as above we obtain that $\F_V$ has one of the following forms:

$$\F_V: \begin{xy}
  \xymatrix{
  [(X, D)] \ \bullet \ar@{<->}[r] & \bullet \ [(X^\vee, D^\vee)] \\
  }
\end{xy} \quad \text{or} \quad \F_V: [(X, D)] \ \bullet \rcirclearrowleft .$$

Moreover, $\F_V$ has a unique vertex if and only if $e = e'$, which is equivalent to $e^2 \equiv 1 \mod d$. As in the proof of Corollary \ref{StructureOfAutomorphismGroup} we see that $\Aut(V_{d, e})$ has the structure of an amalgamated product. To be precise, if $e^2 \not\equiv 1 \mod d$, we let $J := \Aut(V, \pi)$, $J^\vee := \Aut(V, \pi^\vee)$ and $A$ be the subgroup corresponding to the edge. Then, identifying $J \supseteq A \subseteq J^\vee$, we have

$$\Aut(V) = J \star_A J^\vee.$$

In the case $e^2 \equiv 1 \mod d$ we let $A = \langle \Aut(X, D), \psi \rangle$ and $J = \Aut(V, \pi)$, where $\psi: (X, D) \dasharrow (X^\vee, D^\vee) \cong (X, D)$ denotes a reversion. Then

$$\Aut(V) = A \star_{A \cap J} J.$$

This gives an alternative proof of a result of I. Arzhantsev and M. Zaidenberg, see \cite{AZ}, Theorem 4.2.
\end{rem}

\bigskip
\noindent Fakult\"at f\"ur Mathematik, Ruhr-Universit\"at Bochum, Universit\"atstra{\ss}e 150, 44801 Bochum, Germany\\
\textit{E-mail address:} sergei.kovalenko@rub.de 
 
\end{sloppypar}

\end{document}